\documentclass[aos]{imsart}

\RequirePackage{amsthm,amsmath,amsfonts,amssymb}
\RequirePackage[numbers,sort&compress]{natbib}
\RequirePackage[colorlinks,citecolor=blue,urlcolor=blue]{hyperref}
\RequirePackage{graphicx}

\startlocaldefs

\usepackage{subcaption}
\usepackage{float}
\usepackage{verbatim} 

\usepackage{booktabs}
\newcommand{\ra}[1]{\renewcommand{\arraystretch}{#1}}
\usepackage{multirow}
\usepackage{xcolor}

\newtheorem{theorem}{Theorem}[section]
\newtheorem{proposition}[theorem]{Proposition}%
\theoremstyle{thmstyletwo}%
\theoremstyle{thmstylethree}%

\newtheorem{lemma}[theorem]{Lemma}
\newtheorem{corollary}[theorem]{Corollary}

\theoremstyle{definition}
\newtheorem{assumption}{Assumption}

\numberwithin{equation}{section}

\newcommand{\E}{\mathbf{E}}
\newcommand{\D}{D_{\kern-0.2ex F}\kern0.2ex} 
\newcommand{\infdiv}[2]{d_{\kern-0.2ex K\kern-0.2ex L}(#1\kern-0.25ex\left\vert\kern-0.25ex\right\vert#2)}

\newcounter{x}
\def\RN#1{\setcounter{x}{#1}\Roman{x}}
\def\rn#1{\setcounter{x}{#1}\roman{x}}

\newcommand{\vertiii}[1]{{\left\vert\kern-0.25ex\left\vert\kern-0.25ex\left\vert #1 
    \right\vert\kern-0.25ex\right\vert\kern-0.25ex\right\vert}}

\endlocaldefs

\begin{document}

\begin{frontmatter}

\title{Variational non-Bayesian inference of the Probability Density Function in the Wiener Algebra}
\runtitle{Variational non-Bayesian inference of the Probability Density Function}

\begin{aug}
\author[B]{\fnms{U Jin }~\snm{Choi}\ead[label=e2]{ujchoi@kaist.ac.kr}}
\and
\author[A]{\fnms{Kyung Soo}~\snm{Rim}\ead[label=e1]{ksrim@sogang.ac.kr}}

\ \

{October 19, 2023}

\address[B]{Department of mathematical science,
Korea Advanced Institute of Science and Technology \printead[presep={,\ }]{e2}}

\address[A]{Department of mathematics,
Sogang University\printead[presep={,\ }]{e1}}

\end{aug}




\begin{abstract}
This paper presents a research study focused on uncovering the hidden population distribution from the viewpoint of a variational non-Bayesian approach.
It asserts that if the hidden probability density function (PDF) has continuous partial derivatives of at least half the dimension's order, it can be perfectly reconstructed  from a stationary ergodic process:
First, we establish that if the PDF belongs to the Wiener algebra, its canonical ensemble form is uniquely determined through the Fr\'echet differentiation of the Kullback-Leibler divergence, aiming to minimize their cross-entropy.
Second, we utilize the result that the differentiability of the PDF implies its membership in the Wiener algebra.
Third, as the energy function of the canonical ensemble is defined as a series, the problem transforms into finding solutions to the equations of analytic series for the coefficients in the energy function. Naturally, through the use of truncated polynomial series and by demonstrating the convergence of partial sums of the energy function, we ensure the efficiency of approximation with a finite number of data points.
Finally, through numerical experiments, we approximate the PDF from a random sample obtained from a bivariate normal distribution 
and also provide approximations for the mean and covariance from the PDF. This study substantiates the excellence of its results and their practical applicability.
\end{abstract}



\end{frontmatter}

\section{Introduction} \label{intro}

Statistical inference enables us to draw meaningful conclusions about a population's distribution by studying a representative sample.  
Numerous facets of scientific research implicitly or explicitly require the estimation of  population's distribution from a sample. 
By analyzing a subset of individuals or observations, we can gather information that applies to the larger group, even when studying the entire population is impractical or impossible.

Drawing from the works of two authors (\cite{cybenko, hornik}), machine learning methods have gained popularity for predicting hidden information from data.
Today, most emerging artificial neural networks incorporate statistical methods.
Nevertheless, these computations do not achieve global optimization of learning results because of the inherent challenges related to setting initial conditions and handling the small magnitudes of high-dimensional differential values.

To address this issue, many scholars have applied statistical inferences.
Among these methods, variational Bayesian method is a commonly used technique for estimating posterior information.
It's worth noting that there are numerous excellent studies on the method. 
One can refer to the citations \cite{AR,ZG,SBK,BKM} as examples of the most recent works.

Our study aims to identify the hidden probability density function (PDF) without the need for a prior distribution in a non-Bayesian context. This is achieved solely from a stationary ergodic process using entropy minimization. Additionally, we aim to estimate moments of the PDF, such as the mean and variance.
We provide a proof of the norm convergence of sequences of approximated PDFs from a finite sample size, which is useful  in practical applications.
In most probability models arising from natural phenomena, assuming the boundedness of a random variable within a specific  window becomes reasonable when considering a sufficiently wide range of variability.
Moreover, we assume ergodicity in the context of a stationary process, as is common in many scientific communities.

The basic framework of the idea to find the hidden PDF is as follows:
First, we express the PDF in the form of a canonical ensemble, resembling a mechanical heat system.
Second, by embedding the problem in an infinite-dimensional function space, we establish the existence of the Fréchet derivative of an entropy induced from the energy function of the PDF (for more information on the Fr\'echet derivative, refer to \cite{rudin}).
Third, we uniquely determine the energy function by minimizing the Kullback-Leibler divergence (KL-divergence) (\cite{kullback-leibler, csiszar}). The concept of KL-divergence, valuable in information theory, finds widespread use, even in artificial intelligence (for example, see \cite{RHW}, \cite{BT}, \cite{hinton_book2}).
Finally, we derive a system of polynomial series that corresponds to the sample means of complex exponentials and present their numerical solutions
from random samples in a bivariate normal distribution.
Using the approximated PDFs that we obtained, we compare approximations of the mean and variance from the bivariate normal distribution.

\section{Organization and notation}

This article is organized as follows:
In the next section, we formulate the problem by defining the energy and partition function for the hidden PDF in a generalized Wiener algebra. We observe that the nonlinearity of energy functions poses challenges in maximizing likelihood functions. To overcome this issue, we adopt the KL-divergence between two distribution functions.

In Section \ref{sa}, we embed energy functions into a generalized Wiener algebra, thereby transforming the nonlinearity of the energy function into a system of equations in the function space and simultaneously converting the min-max problem into a problem of solving equations.
We show the Fr\'echet differentiability of the KL-divergence with respect to coefficients in the generalized Wiener algebra. From this property, we then uniquely obtain the existence such that its energy function yields a KL-divergence of $0$.

We establish a characterization of the coefficients of the energy function by expressing them in relation to the Fourier coefficients found within the (classical) Wiener algebra, as demonstrated in Section \ref{pnn}. This serves as an example of the generalized Wiener algebra.

In the concluding Section \ref{main section}, when dealing with a random sample that consists of realizations of a stationary ergodic process, we formulate a system of polynomial series with the coefficients of the energy function as variables, which consists of infinitely many equations in infinitely many variables.
From the convergence property of the partial sums of an energy function in Section \ref{sa}, the truncated system, consisting of a finite number of equations with partial sums, offers an avenue for coefficient approximations.
The estimated coefficients effectively serve as estimators for the coefficients. More precisely, the hidden PDF can be approximated using these estimators.
Furthermore, we provide a numerical example that emphasizes the approximation of the PDF, as well as the calculation of the mean and variance from the PDF using a random sample. This sample is generated from a bivariate normal distribution.

In this paper, we consistently employ the following general notations (Table \ref{T0}). Some notations will be elaborated upon in more detail when they are introduced.

\begin{table}[h]\centering
\ra{1} 
\begin{tabular}{ll}\toprule
$\mathbb{R}^d$	&  the $d$-dimensional Euclidean space \\ 
$\mathbb{Z}^d$	&  the additive group of all $d$-tuples of integers, which \\ & serves as a multi-index set  \\ 
$|\alpha|_\infty$	&   $\max_{1\le k\le d} |\alpha_k|$ for $\alpha=(\alpha_1,\ldots,\alpha_d)\in\mathbb{Z}^d$ \\
$S$	& a compact subset of $\mathbb{R}^d$ \\
$\mathbb{T}^d$             &  the $d$-torus $[-\pi,\pi]^d$ is equipped with the quotient topology \\ & obtained by identifying opposite edges \\ 
$X$ 		& $X=(X_1,\ldots,X_d)$ is a random vector   \\ 
$X_k$             &   the $k$-th component of $X$ (a random variable) \\
$X[n]$     &  the random vector $X[n]=(X_1[n],\ldots,X_d[n])$ in a sequence    \\
$X_k[n]$     &  the $k$-th component of $X[n]$  (a random variable) \\
$x$ 		& $x=(x_1,\ldots,x_d)$ is a realization of $X$   \\ 
$x_k$             &   the $k$-th scalar component of $x$ \\
$x[n]$     &  $x[n]=(x_1[n],\ldots,x_d[n]) $ is a realization of $X[n]$ in a sequence   \\
$x_k[n]$     &  the $k$-th scalar component of $x[n]$ \\
$p_0(x)$   &  the hidden PDF of $X$    \\
$C(S)$   &   the set of all continuous functions on $S$    \\
$\ell^1$   &  the space of summable complex sequences on $\mathbb{Z}^d$  \\
			& with norm $\|\cdot\|_{\ell^1}$ \\
$c_0(\mathbb{C})$   &  the subspace of $\ell^\infty(\mathbb{C})$, where each element \\ & possesses a finite number of supports  \\
$y$   &  a sequence of complex numbers $y=(y_\alpha)\in \ell^1(\mathbb{C})$ \\ & which  serves as a coefficient sequence of $E(x,y)$\\ 
$E(x,y)$   &  the energy function of $x$ with a coefficient sequence $y$ \\
$p(x\mid y)$   &  a PDF of $X$ with $y$ ({\it not} a conditional PDF) \\
$Z_E$   &  the partition function of $E$   \\
$\mathbf{E}_p[X]$   &  the expectation of $X$ which has the PDF $p$  \\
$\infdiv{p_1}{p_2}$   &  the Kullback-Leibler divergence between the PDFs $p_1$ and $p_2$    \\
$H(p_1,p_2)$   &  the cross-entropy of $p_2$ relative to $p_1$ \\ 
$H(p)$   &   the entropy of $p$    \\
$\mathcal{A}(\phi_\alpha)$   &   a generalized Wiener algebra   \\
$dx$   &   the standard Lebesgue measure on $\mathbb{R}^d$    \\
$dm_d$   &   the normalized Lebesgue measure $dm_d=dx/(2\pi)^d$    \\
$\hat f(\alpha)$   &   the $\alpha$-th Fourier coefficient of $f$, where $\alpha\in\mathbb{Z}^d$    \\
$E(\mathbf{0},y)^{\ast n}$   &  the $n$-fold iteration of convolution with itself \\ & ($\mathbf{0}\in\mathbb{R}^d$, $y\in\ell^r(\mathbb{C})$)   \\
$\delta_{\mathbf{0}}$   &   the Dirac delta concentrated at $\mathbf{0}$  (cf. $\delta_{0}$  at $0\in\mathbb{R}$)  \\
$\delta_{\alpha\beta}$   &   the Kronecker delta for multi-indices $\alpha$, $\beta$\\
$\left<\,\cdot\,,\,\cdot\,\right>_{L^2(S)}$   &   the inner product in $L^2(S)$  \\
$\lceil\cdot\rceil$ & the ceiling function \\
$:=$   &   the expression on the left-hand side is defined to be equal to \\ & the right-hand side \\\bottomrule
\end{tabular}
\caption{Notations}
\label{T0}
\end{table}

\section{Problem Setting and Preliminaries} \label{epf}

This section will set the mathematical model that is being addressed and provide basic background information relevant to understanding it.

Initially, let's assume that the random vector $X = (X_1, ..., X_d)$ has a PDF $p_0$ and is confined within a suitable window, as explained in Section \ref{intro}.
To clarify, each observed vector $x = (x_1, ..., x_d)$ exists within a window of a compact set $S \subset \mathbb{R}^d$. 
In this article, the information regarding $p_0$ is completely unknown.

Next, we define an auxiliary PDF $p(x\mid y)$ with a countable coefficient set $y=(y_\alpha)$ by
\begin{equation} \label{network probability-1}
     p(x\mid y) = \frac{e^{-E(x,y)}}{Z_E},
\end{equation}
where, $E(x,y)$ is referred to as the energy function, and $Z_E=\int_S e^{-E}dm_d$ is the partition function for $p:=p(x\mid y)$, which takes the form of a canonical ensemble, similar to the Boltzmann distribution. The auxiliary PDF $p$ is entirely determined by $E$.
The objective of this study is to find a sequence of coefficients $y$ that satisfies the following equation  based on the differentiability condition of $p_0$ (refer to Assumption \ref{A} below):
\begin{equation} \label{goal}
p(x\mid y) = p_0(x).
\end{equation}

If all components $X_1, \ldots, X_d$ of $X$ are independent random variables from $p_0$. Assuming  (\ref{goal}), the joint PDF of $X_1, \ldots, X_d$ is equal to the product of the individual PDFs of $X_1,\ldots,X_d$. 
Thus, $E(x, y)$ can be expressed as
\begin{equation} \label{iid}
E(x, y) = E(x_1, y^{(1)}) + \cdots + E(x_d, y^{(d)})
\end{equation}
for some $y^{(k)}:=(y^{(k)}_\eta)$ $(\eta\in\mathbb{Z})$ $(1\le k\le d)$. 
Under the independent and identically distributed (i.i.d.) condition, we have $y^{(1)}=\cdots=y^{(d)}$.

Throughout this study, the following assumption is made:

\begin{assumption} \label{A}
	All of $p_0$'s partial derivatives of order $\lceil d/2 \rceil$ are continuous, and the entropy of $p_0$ is finite, i.e., $|H(p_0)| < \infty$.
\end{assumption}
When performing logarithmic calculations of $p_0$, the points where $p_0$ becomes zero might cause discontinuities.
To circumvent this, since $S$ is compact, we will solve the equation of  $p = p_0 + 1$ instead of $p=p_0$ in (\ref{goal}). 
Subsequently, we infer $p_0$ by substituting $p - 1$ for $p$.
Thus, we reaffirm Assumption \ref{A'} instead of Assumption \ref{A} solely in the strict context of mathematical analysis:

\begin{assumption} \label{A'}
On a window, $p_0 > 0$,  all of its partial derivatives of order $\lceil d/2 \rceil$ are continuous, and  the entropy of $p_0$ is finite, i.e., $|H(p_0)| < \infty$.
\end{assumption}

Now, let's introduce the elementary properties of KL-divergence. For two PDFs $p_1$ and $p_2$, their KL-divergence is defined as
\begin{equation*}
\infdiv{p_1}{p_2} = \E_{p_1}\left[\ln\frac{p_1}{p_2}\right],
\end{equation*}
which is defined for $x$ such that $p_1(x)p_2(x)\ne0$.
The KL-divergence satisfies the following three conditions:

\begin{enumerate}
\item[$(\rn{1})$] $\infdiv{p_1}{p_2}\ge0$ \quad (Gibbs' inequality).
\item[$(\rn{2})$] $\infdiv{p_1}{p_2}=0$\, if and only if \,$p_1=p_2$\; almost everywhere \quad (identity of indiscernible).
\item[$(\rn{3})$] $\infdiv{p_1}{p_2}\ne \infdiv{p_2}{p_1}$ \quad (asymmetricity).
\end{enumerate}
Note that, more generally, for non-negative measurable functions $f_1$ and $f_2$, the equality condition of Jensen's inequality can be extended to $\infdiv{f_1}{f_2}=0$ if and only if
\begin{equation} \label{extended identity}
f_1=cf_2 \quad \text{almost everywhere}
\end{equation}
for some constant $c$.

The KL-divergence for the two mentioned PDFs at the beginning of this section is as follows:
\begin{equation} \label{cross entropy}
\begin{aligned}
    \infdiv{p_0}{p}
        &=-\E_{p_0}[\ln p(x\mid y)] - (-\E_{p_0}[\ln p_0]) \\ 
        &= H(p_0,p) - H(p_0)\ge0,
\end{aligned}
\end{equation} 
where $H(p_0,p)$ and $H(p_0)$ denote the cross-entropy of $p$ relative to $p_0$, and the entropy of $p_0$, respectively. Consequently, for all $y$, $H(p_0,p)\ge H(p_0)$ holds from Assumption \ref{A'} . 

To solve (\ref{goal}), using property $(\rn{2})$, we only need to find $y$ such that $\infdiv{p_0}{p}=0$.
Since the quantities $\infdiv{p_0}{p}$ and $H(p_0,p)$ are equal up to the constant difference $H(p_0)$, the first step towards finding the most efficient solution is to determine $y$ that minimizes
\begin{equation} \label{entropy form}
     \infdiv{p_0}{p},  
\end{equation}
which is equal to minimize
\begin{equation} \label{entropy form2}
     H(p_0,p),
\end{equation}
where two quantities conceive $y$, since $p$ depends on both $x$ and $y$.

Unfortunately, directly clarifying (\ref{entropy form}) or (\ref{entropy form2}) is difficult due to the nonlinearity of $H(p_0,p(x\mid y))$ generally. (In \cite{barron-sheu}, a decaying estimate of $\infdiv{p_0}{p}$ is obtained.) 
To address this issue, in the next section, we will expand a nonlinear energy function within a generalized Wiener algebra, which is a type of Wiener algebra.

\section{Generalized Wiener algebra and Fr\'echet derivative} \label{sa}

In this section, we define a space of energy functions and demonstrate that $\infdiv{p_0}{p}$ is Fr\'echet differentiable 
with respect to each $y_\alpha$ on it.

Suppose that $(\phi_\alpha)$ is a collection of orthonormal continuous functions, i.e., $\left<\phi_\alpha,\phi_\beta\right>_{L^2(S)}=\delta_{\alpha\beta}$, such that 
\begin{equation} \label{group}
	\phi_\alpha\phi_\beta=\phi_{\alpha+\beta}
\end{equation}
for all multi-indices and such that
\begin{equation} \label{character}
	\vertiii{(\phi_\alpha)}_{\infty}
		:=\sup_{x\in S}\sup_{\alpha\in\mathbb{Z}^d}|\phi_\alpha(x)|<\infty.
\end{equation}

Now we define a class of functions by 
\begin{equation*} 
	\mathcal{A}(\phi_\alpha)=\left\{ \sum_{\alpha\in\mathbb{Z}^d} y_\alpha\phi_\alpha(x) 
				\,\Big\vert\, 
				(y_\alpha)\in\ell^1\right\},
\end{equation*}
where $\ell^1$ is the space of 
summable sequences $(y_\alpha)$, i.e.,  
$\|(y_\alpha)\|_{\ell^1}=\sum_{\alpha\in\mathbb{Z}^d} |y_\alpha| < \infty$.
For further information regarding the $\ell^1$-space, refer to \cite{rudin}, for instance.

We briefly mention the properties of $\mathcal{A}(\phi_\alpha)$ that are required for the proof in the main theorems:
By the triangle inequality and by H\"older's inequality, we have $|E(x,y)|\le\sum_{\alpha\in\mathbb{Z}^d}|y_\alpha\phi_\alpha|$, which is bounded  by $\vertiii{(\phi_\alpha)}_{\infty}\|(y_\alpha)\|_{\ell^1}<\infty$. Therefore, the sum converges pointwise.
Moreover, 
\begin{equation} \label{uniform convergence property}
\left|\sum_{\alpha\in\mathbb{Z}^d}y_\alpha\phi_\alpha-\sum_{|\alpha|_\infty\le N}y_\alpha\phi_\alpha \right|
	\le \left|\sum_{|\alpha|_\infty> N}y_\alpha\phi_\alpha \right| 
	\le \vertiii{\phi_\alpha}_{\infty}\sum_{|\alpha|_\infty> N}|y_\alpha| \to0
\end{equation}	
as $N\to\infty$. That is, the partial sums $\sum_{|\alpha|_\infty\le N}y_\alpha\phi_\alpha$ of $E$ converges uniformly on $S$, and furthermore, they converge absolutely as well.

By applying the inner product, we can establish that the series representation of $E$ in $\mathcal{A}(\phi_\alpha)$ is unique.
Indeed, if there are two representations, namely $E := \sum_\alpha y_\alpha\phi_\alpha = \sum_\alpha y'_\alpha\phi_\alpha:=E'$, then $y_\alpha=\left<E,\phi_\alpha\right>_{L^2(S)}=\left<E',\phi_\alpha\right>_{L^2(S)}=y'_\alpha$, which holds for all $\alpha$ by  the Lebesgue dominated convergence theorem via (\ref{character}) and (\ref{uniform convergence property}).

The class $\mathcal{A}(\phi_\alpha)$, equipped with the norm $\|E\|_{\mathcal{A}} := \|y\|_{\ell^1}$, 
forms a Banach space due to the completeness of $\ell^1$. 
This space is isometrically isomorphic to $\ell^1$ and has a pre-dual $c_0$ consisting of vanishing sequences. (It's worth noting that the space $c_c$ of sequences with only finitely many nonzero elements is dense in both $\ell^1$ and $c_0$.)

Using (\ref{group}) and applying Fubini's theorem for infinite series, we find that
\begin{equation*}
\begin{split}
	E(x,y)E'(x,y') 
	&= \sum_{\alpha\in\mathbb{Z}^d}\sum_{\beta\in\mathbb{Z}^d}y_\alpha y'_\beta \phi_{\alpha+\beta}(x) \\
	&= \sum_{\beta\in\mathbb{Z}^d}\sum_{\alpha\in\mathbb{Z}^d}y_\alpha y'_{\beta-\alpha} \phi_{\beta}(x) \\
	&= \sum_{\beta\in\mathbb{Z}^d} (y\ast y')(\beta)\phi_{\beta}(x). 
\end{split}
\end{equation*}
Furthermore, as $\|y\ast y'\|_{\ell^1} \le \|y\|_{\ell^1} \|y'\|_{\ell^1}$, we deduce $\|EE'\|_{\mathcal{A}} \le \|E\|_{\mathcal{A}} \|E'\|_{\mathcal{A}}$, solidifying its status as a Banach algebra. This serves as a generalization of the Wiener algebra.

If $E^{(n)}$ is a bounded sequence in $\mathcal{A}(\phi_\alpha)$ and converges to $E$ uniformly on $S$, then employing the Banach-Alaoglu theorem allows us to identify a subsequence $E^{(n_k)}$ that converges in the weak-$\ast$ sense with a duality pairing  $\left<\,\cdot\,,\,\cdot\,\right>$ between $\ell^p$-spaces. This subsequence, alongside an element $E' \in \mathcal{A}^r(\phi_\alpha)$, satisfies $\lim_{k \to \infty} \left<h, E^{(n_k)}\right> = \left<h, E'\right>$ for any $h\in c_0$.
By choosing $h$ such that only $h_\alpha=1$ and $0$ otherwise, it can be inferred that the $\alpha$-th coefficient of $E^{(n_k)}$ converges to the corresponding coefficient of $E'$.
But since $E^{(n_k)}$ converges uniformly to $E$, the coefficients of $E^{(n_k)}$ also converge pointwise to the coefficients of $E$.
By the unique representation property above, we conclude $E=E'\in \mathcal{A}(\phi_\alpha)$.

Therefore, the Banach algebra $\mathcal{A}(\phi_\alpha)$ is summarized as follows:

\begin{proposition}\label{generalized wiener algebra} Suppose that $E \in \mathcal{A}(\phi_\alpha)$. Then
\begin{enumerate}
\item[$(a)$]  $E$ has a uniform upper bound, 
	$|E(x,y)|\le \vertiii{(\phi_\alpha)}_{\infty}\|(y_\alpha)\|_{\ell^1}<\infty$.
\item[$(b)$] The partial sums of $E$ converge uniformly and absolutely on $S$.
\item[$(c)$] The representation of $E$ is unique.
\item[$(d)$] If $E^{(n)}$ is a bounded sequence in $\mathcal{A}(\phi_\alpha)$ and converges to $E$ uniformly on $S$, then $E\in \mathcal{A}(\phi_\alpha)$.
\end{enumerate}
\end{proposition}

To state  the main theorem of this section we need the following lemma.

\begin{lemma} \label{convexity}
If $E\in\mathcal{A}(\phi_\alpha)$, then $\infdiv{p_0}{p}$ is Fréchet differentiable on $\ell^1$. Moreover, its Fréchet derivative is given by
\begin{equation} \label{frechet}
 \D\infdiv{p_0}{p}(h) = \sum_{\alpha\in\mathbb{Z}^d} \frac{\partial}{\partial{y_\alpha}} H(p_0,p) h_\alpha
\end{equation}
for $h\in\ell^1$.
\end{lemma}

\begin{proof}
Let $y=(y_\alpha)\in\ell^1$. Considering (\ref{cross entropy}), we focus on $H(p_0,p)$ instead of $\infdiv{p_0}{p}$ since $H(p_0)$ is a finite number. This leads to the following expansion of the cross entropy:
\begin{equation} \label{expansion of cross entropy}
\begin{aligned}
	H(p_0,p) 
	&=\E_{p_0}[E(x,y) + \ln Z_E] \\
	&=\E_{p_0}\Big[\sum_{\alpha\in\mathbb{Z}^d}  y_\alpha\phi_\alpha(x)\Big] 
		+ \ln \int_{S} \!e^{-\sum_{\alpha\in\mathbb{Z}^d}  y_\alpha\phi_\alpha(x)} dx.	
\end{aligned}
\end{equation}
Here, first note that $H(p_0,p)$ is well-defined, implying finiteness for all $y$. Due to the compactness of $S$ via $(a)$ of Proposition \ref{generalized wiener algebra}, 
two terms of (\ref{expansion of cross entropy}) are readily finite. 

We will show that the $y_\alpha$-partial derivative of $H(p_0,p)$ exists. By interchangeability of the integral and differentiation in (\ref{expansion of cross entropy}), we have
\begin{equation} \label{F-deriv}
\begin{aligned} 
	\frac{\partial}{\partial{y_\alpha}}H(p_0,p) 
	        =\E_{p_0}\left[ \frac{\partial}{\partial{y_\alpha}} E(x,y) \right]
            		+ \frac{\frac{\partial}{\partial{y_\alpha}}Z_E}{Z_E}, 
\end{aligned}	 
\end{equation}	
which is equal to
\begin{equation} \label{f-derivative}
\begin{aligned} 
	        &\E_{p_0}\left[\frac{\partial}{\partial{y_\alpha}} E(x,y) \right]
			- \frac{1}{Z_E}\int_S e^{-E(x,y)}\frac{\partial}{\partial{y_\alpha}} E(x,y)\,dx  \\
	        &=\E_{p_0}\left[\frac{\partial}{\partial{y_\alpha}} E(x,y) \right]
			- \E_{p}\left[\frac{\partial}{\partial{y_\alpha}} E(x,y) \right] \\
	        &=\E_{p_0}\left[\phi_\alpha \right]
			- \E_{p}\left[\phi_\alpha \right].
\end{aligned}	 
\end{equation}	
In addition, the last terms in (\ref{f-derivative}) are finite. Indeed, those are bounded by $2\|\phi_\alpha\|_{L^\infty}\le 2\vertiii{(\phi_\alpha)}_{\infty}<\infty$ for every $\alpha$. 
Thus, $\frac{\partial}{\partial{y_\alpha}} H(p_0,p)$ is derived for any $\alpha$.

Now we will prove that the Fr\'echet derivative of $H(p_0,p)$ is induced by
\begin{equation*}
	\D H(p_0,p)(h) = \sum_{\alpha\in\mathbb{Z}^d}  \frac{\partial}{\partial{y_\alpha}} H(p_0,p)h_\alpha
\end{equation*} 
as a bounded linear functional of $h\in\ell^1$.
Since $c_c$ is a dense subset of  $\ell^1$, thus, we use $c_c$ instead of  $\ell^1$.

Let $h\in c_c$ and put
\begin{equation*} 
	H(p_0,p(x\mid y+h))- H(p_0,p(x\mid y)) 
			- \sum_{\alpha\in\mathbb{Z}^d} \frac{\partial}{\partial{y_\alpha}} H(p_0,p)  h_\alpha:=\RN{1}.
\end{equation*}
First, we will show that $\RN{1}=o(\|h\|_{\ell^1})$ as $\|h\|_{\ell^1}\to0$. 
It follows that
\begin{equation*} 
	\RN{1} = \sum_{\alpha\in\mathbb{Z}^d}   
			\left(\int_0^1 \frac{\partial}{\partial{y_\alpha}}  H(p_0,p(x\mid y+th)) \,dt - \frac{\partial}{\partial{y_\alpha}} H(p_0,p) \!\right) h_\alpha,	
\end{equation*}
where the equality is due to the chain rule as the summation runs over only a finite number of $\alpha$.

In substituting (\ref{F-deriv}) via (\ref{f-derivative}) into $\RN{1}$, by the triangle inequality, we obtain
\begin{equation} \label{continuity-vanishing}
\begin{split}
	|\RN{1}| 
	&\le \sum_{\alpha\in\mathbb{Z}^d}\int_0^1\!\! \int_S |h_\alpha\phi_\alpha(x) | | p(x\mid y+th) -  p(x\mid y) | dx dt \\
	&\le \vertiii{(\phi_\alpha)}_\infty \|h\|_{\ell^1} \int_0^1\!\! \int_S  |p(x\mid y+th) - p(x\mid y)| dx dt,
\end{split}		
\end{equation}
where the second inequality follows from H\"older's inequality. 
By applying the multi-dimensional mean value theorem for integrals, we obtain
\begin{equation} \label{multi-mvt-pre}
\begin{split}
	\int_0^1\int_S  & |p(x\mid y+th) - p(x\mid y)| dx dt \\
		&= \int_0^1 t \int_S  
			|\nabla_{y} p(x\mid y+\xi_{th})\cdot h| dx dt, 
\end{split}		
\end{equation}
where $\xi_{th}$ denotes some vector in the ball with radius $\|h\|_{\ell^2}$ (the ball  is finite dimensional in a Euclidean space) 
and $\nabla_y$ represents the ordinary gradient for those specific indices of $y_\alpha$ where $h_\alpha \neq 0$.
By H\"older's inequality, the last integral (\ref{multi-mvt-pre}) is bounded by
\begin{align}
	\|h\|_{\ell^1} \int_0^1 t \int_S  
	& \|\nabla_{y}p(x\mid y+\xi_{th})\|_{\ell^{\infty}}dx dt  \notag \\
	&= \|h\|_{\ell^1} \int_0^1 t \int_S  
		\sup_{\alpha\in \{\alpha | h_\alpha\ne0\} } \Big|\frac{\partial}{\partial{y_\alpha}} p(x\mid y+\xi_{th})\Big| dx dt. \label{multi-mvt}
\end{align}		

Now we will estimate (\ref{multi-mvt}). For simplicity, put $Z=Z_{E(x,y+\xi_{th})}$.
By interchangeability of the integral and differentiation, 
\begin{equation*}
\begin{split}
\frac{\partial}{\partial y_\alpha} p(x\mid y+\xi_{th}) 
	&= \frac{\partial}{\partial y_\alpha} \frac{e^{-E(x,y+\xi_{th})}}{Z} \\
	&= -\phi_\alpha(x)p(x\mid y+\xi_{th})
		- p(x\mid y+\xi_{th})\int_S\frac1{Z} \frac{\partial}{\partial y_\alpha} e^{-E(x,y+\xi_{th})}dx \\
	&= -\phi_\alpha(x)p(x\mid y+\xi_{th})
		+ p(x\mid y+\xi_{th})\int_S\phi_\alpha(x)p(x\mid y+\xi_{th})dx 
\end{split}
\end{equation*}
which is bounded by $2\vertiii{(\phi_\alpha)}_\infty p(x\mid y+\xi_{th})$.
Thus, the iterated integral of  (\ref{multi-mvt}) is bounded by $\vertiii{(\phi_\alpha)}_\infty$.

According to (\ref{continuity-vanishing}) $\sim$ (\ref{multi-mvt}), we get $| \RN{1} | \le \vertiii{(\phi_\alpha)}_\infty^2\|h\|_{\ell^1}^2$. Consequently, $| \RN{1} | / \|h\|_{\ell^1} \to0$ as $\|h\|_{\ell^{r}}\to0$.
Hence,
\begin{equation} \label{functional}
	\D H(p_0,p)(h) = \sum_{\alpha\in\mathbb{Z}^d} \frac{\partial}{\partial{y_\alpha}} H(p_0,p) h_{\alpha}
\end{equation}
holds for any $h\in c_c$.

Finally, it suffices to prove that $\D H(p_0,p)$ is a bounded linear functional on $c_c$ with norms $\|\cdot\|_{\ell^1}$. 
%
Choose $h\in c_c$ such that $\|h\|_{\ell^r}=1$.   Combining (\ref{functional}), (\ref{F-deriv}), and (\ref{f-derivative}),  we have
\begin{equation*}
\begin{aligned}
	|\D H(p_0,p)(h)| 
	&\le \left\|  \left(\frac{\partial}{\partial{y_\alpha}}H(p_0,p)\right)\right\|_{\ell^\infty} \|h\|_{\ell^1} \\
	&= \sup_{\alpha\in\mathbb{Z}^d} \big|\E_{p_0}[\phi_\alpha] - \E_p[\phi_\alpha] \big| \\
	&\le \vertiii{(\phi_\alpha)}_\infty \int_S |p_0(x)-p(x\mid y)|dx \\
	&\le2\vertiii{(\phi_\alpha)}_\infty <\infty.
\end{aligned}
\end{equation*}
This completes the proof.
\end{proof}

In the proof of Lemma \ref{convexity}, the condition that (\ref{f-derivative}) vanishes for some $y$, i.e., $\D H(p_0,p(x\mid y))=0$, 
is a necessary condition for $y$ to be a solution to (\ref{goal}).
If $\ln p_0$ belongs to the generalized Wiener algebra, the following theorem shows the unique existence of such $y$.

\begin{theorem} \label{basic equations}
If $\ln p_0\in\mathcal{A}(\phi_\alpha)$, then
there exists a unique solution $y\in\ell^1$ to the equations of 
$\E_{p_0}[\phi_\alpha] = \E_{p(x\mid y)}[\phi_\alpha]$ for every $\alpha$, at which $(\ref{goal})$ holds.
Moreover, $H(p_0) = H(p(x\mid y))$.
\end{theorem}

\begin{proof} 
Let $E\in\mathcal{A}(\phi_\alpha)$ and put $p(x\mid y)=e^{-E(x,y)}/Z_E$. 
First, we will show $e^{-E}\in \mathcal{A}(\phi_\alpha)$. 
Since $\mathcal{A}(\phi_\alpha)$ is a Banach algebra, the partial sums of $\sum_{k=0}^n \frac{(-1)^k}{k!}E^k$ belongs to $\mathcal{A}(\phi_\alpha)$.
By the property $(a)$ of Proposition \ref{generalized wiener algebra} we have $\|E^k\|_{\mathcal{A}}\le \|E\|^k_{\mathcal{A}}\le\vertiii{(\phi_\alpha)}^k_\infty\|y\|^k_{\ell^1}$. Thus,
\begin{equation} \label{uniform exponential}
	\Big\|\sum_{k=0}^n\frac{(-1)^k}{k!}E^k\Big\|_{\mathcal{A}}
	\le \sum_{k=0}^n\frac{1}{k!}\|E\|^k_{\mathcal{A}} 
	\le e^{\vertiii{(\phi_\alpha)}_\infty \|(y_\alpha)\|_{\ell^1}}.
\end{equation}
Consequently,  $\sum_{k=0}^n\frac{(-1)^k}{k!}E^k$ is
a bounded sequence in $\mathcal{A}(\phi_\alpha)$. 
Since the tail of $e^{-E} = \sum_{k=0}^\infty \frac{(-1)^k}{k!}E^k$ vanishes uniformly on $S$, 
according to property $(d)$ of Proposition \ref{generalized wiener algebra}, $e^{-E}$ is an element of $\mathcal{A}(\phi_\alpha)$. So is
$p(x\mid y)$. 
Similarly, we can justify that $p_0\in\mathcal{\phi_\alpha}$ from the assumption that $\ln p_0\in \mathcal{A}(\phi_\alpha)$.

Let $y\in\ell^1$ be a solution to $\E_{p_0}[\phi_\alpha] - \E_{p(x\mid y)}[\phi_\alpha]=\int_S (p_0(x)-p(x\mid y))\phi_\alpha(x)dx=0$ for every $\alpha$. 
Also, since $p_0(x)-p(x\mid y)\in\mathcal{A}(\phi_\alpha)$, it has a unique representation by $(c)$ of Proposition \ref{generalized wiener algebra}.
Hence, $p(x\mid y)-p_0(x)=0$ on $S$, which leads to $(\ref{goal})$. 
 
To demonstrate the uniqueness of $y$, let's assume that $\E_{p(x\mid y)}[\phi_\alpha] = \E_{p(x\mid y')}[\phi_\alpha]$  for every $\alpha$. 
Employing a similar approach as before, we can deduce
$p(x\mid y) = p(x\mid y')$ and subsequently $\ln p(x\mid y) = \ln p(x\mid y')$ throughout $S$, given that the logarithm is strictly monotonic.
Thus, their energy functions $E(x,y)$ and $E(x,y')$ are identical within $\mathcal{A}(\phi_\alpha)$. 
Finally, according to $(c)$ of Proposition \ref{generalized wiener algebra}, we have $y_\alpha = y'_\alpha$ for every $\alpha$.


Therefore, the proof is now complete.
\end{proof}

Theorem \ref{basic equations} establishes the unique existence of coefficients $y_\alpha$ satisfying the objective (\ref{goal}) of the paper. However, 
the hypothesis of the theorem requires that, in the absence of knowledge about $p_0$, $\ln p_0$ must be included in the generalized Wiener algebra.
Furthermore, the theorem does not offer a  concrete way to determining the coefficients $y_\alpha$.

Nevertheless, Theorem \ref{basic equations} conveys its significance by providing guidance for selecting the desired $\phi_\alpha$ such that 
$\ln p_0\in\mathcal{A}(\phi_\alpha)$. More precisely, if we can identify a Banach algebra to which $\ln p_0$ belongs and find a suitable Schauder basis $\phi_\alpha$ that generates this space, it's possible that we could eliminate the hypothesis of Theorem \ref{basic equations}. 
In the upcoming sections, we will introduce the Wiener algebra and aim to eliminate the hypothesis based on Assumption \ref{A'} as stated in Section \ref{epf}. We will also cover the process of concretizing and approximating the coefficients $y_\alpha$ while showcasing a method for determining them.

\section{Characterization of coefficients} \label{pnn} 

In this section, we consider the compact set $S = \mathbb{T}^d = [-\pi, \pi]^d$ to be a $d$-dimensional torus, serving as a window for random vectors. 
As a typical example of energy functions, we can consider the Wiener algebra of Fourier series, which is presented as
$\mathcal{A}(e^{i \alpha\cdot x})$.

We specifically utilize the Wiener algebra to elucidate the nature of $y_\alpha$ and establish the convergence of partial sums for an energy function.
More precisely,
if $f\in \mathcal{A}(e^{i \alpha\cdot x}) $, then $f(x)=\sum_{\alpha}\hat{f}(\alpha) e^{i\alpha\cdot x}$ for some $(\hat f(\alpha))\in\ell^1$, where
$\hat{f}(\alpha)=\int_{\mathbb{T}^d} f(x)e^{-i\alpha\cdot x}dm_d $ represents the the $\alpha$-th Fourier coefficient of $f$ and $dm_d=dx/(2\pi)^d$. 
Note that $\vertiii{(e^{i \alpha\cdot x})}_{\infty}=1<\infty$.
For further details about the Wiener algebra, refer to the references \cite{KATZ} and \cite{wiener}.

The proposed energy function in this study is expressed as an infinite series, making it impractical to compute all $y_\alpha$. However, the following result provides an approximation of the auxiliary PDF constructed from partial sums of the energy function, thereby proving its $L^1$-convergence.

\begin{theorem} \label{decaying}
Let $E_N=\sum_{|\alpha|_\infty\le N} y_\alpha e^{i \alpha\cdot x}$ be the $N$-th partial sum of $E$ in $\mathcal{A}(e^{i \alpha\cdot x})$. 
Then
\begin{equation*}
		\lim_{N\to\infty}\int_{\mathbb{T}^d} \left| \frac{e^{-E(x,y)}}{Z_E} - \frac{e^{-E_N(x,y)}}{Z_{E_N}}\right|dx = 0.
\end{equation*}
%
\end{theorem}

\begin{proof}
By Taylor's series expansion of the exponential function, the triangle inequality, and the uniform convergence, we have
\begin{equation} \label{L^1-convergence}
\begin{aligned}
	&\int_{\mathbb{T}^d} \left|e^{-E(x,y)-\ln Z_E} - e^{-E_N(x,y)-\ln Z_{E_N}} \right| dx \\
	&\le \sum_{n=1}^\infty\frac{1}{n!}\int_{\mathbb{T}^d} e^{-E(x,y) - \ln Z_E } | E(x,y) - E_N(x,y) + \ln Z_E/Z_{E_N} |^n dx \\
	&\le \sum_{n=1}^\infty\frac{1}{n!}\left(\sum_{|\alpha|_\infty>N} |y_\alpha| + |\ln Z_E/Z_{E_N} |\right)^n\int_{\mathbb{T}^d} e^{-E(x,y)-\ln Z_E}dx \\
	&=\sum_{n=1}^\infty\frac{1}{n!}\left(\sum_{|\alpha|_\infty>N} |y_\alpha|+ |\ln Z_E/Z_{E_N} |\right)^n \\
	&= e^{\sum_{|\alpha|_\infty>N} |y_\alpha|+ |\ln Z_E/Z_{E_N} |}-1.
\end{aligned}	
\end{equation}

In (\ref{L^1-convergence}), since $\sum_{|\alpha|_\infty>N} |y_\alpha|$ goes to 0 as $N\to\infty$, the term
\begin{equation*}
	Z_{E_N} = \int_{\mathbb{T}^d} e^{-E_N}dx \to Z_E \quad\mbox{as } N\to\infty.
\end{equation*}
Therefore, the last term of (\ref{L^1-convergence}) vanishes as $N\to\infty$, and the proof is complete.
\end{proof}

We will now proceed to show the main result and provide its proof in this section.

\begin{theorem}\label{PNN by frequency}
Assuming Assumption \ref{A'}, $y_\alpha$ is uniquely determined by 
\begin{equation} \label{infor} 
 	y_\alpha = \widehat{\ln \frac1{p_0}}(\alpha) 
\end{equation}
in order to satisfy $(\ref{goal})$, where $\widehat{\ln \frac1{p_0}}(\alpha)$ represents the $\alpha$-th Fourier coefficient of $\ln \frac1{p_0}$.
\end{theorem}

\begin{proof}
Based on the premise denoted as Assumption \ref{A'}, $p_0$ has continuous partial derivatives of order 
$\lceil d/2\rceil$. By Theorem 3.2.16 of \cite{grafakos}, $p_0$ possesses a unique Fourier series expansion characterized by uniform convergence, this property can also be extended to $\ln p_0$, since $p_0>0$. That is, $\ln p_0$ belongs to $\mathcal{A}(e^{i \alpha\cdot x})$. According to Theorem \ref{basic equations}, there exists a unique $y$ such that $\E_{p_0}[\phi_\alpha] = \E_{p(x\mid y)}[\phi_\alpha]$ for every $\alpha$, ensuring that $(\ref{goal})$ is satisfied. Therefore, we have $y_\alpha = \widehat{\ln\frac1{p}}(\alpha) = \widehat{\ln\frac1{p_0}}(\alpha)$, which completes the proof.
\end{proof}

If all components $X_k$ of $X$ are independent, 
then the joint PDF of $X$ becomes $p_0(x)=\prod_{k=1}^dp_{0,k}(x_k)$, where $p_{0,k}$ is the marginal PDF of $X_k$.
By Theorem \ref{PNN by frequency}, we have
\begin{equation} \label{iid fourier coefficient}
\begin{split}
	y_\alpha 
		&= \sum_{k=1}^d \int_{\mathbb{T}^d} e^{-i\alpha\cdot x}\ln \frac{1}{p_{0,k}(x_k)}dm_d \\
		&= \left\{
		      \begin{array}{ll}
 		       \sum_{k=1}^d \int_{\mathbb{T}^1} \ln \frac{1}{p_{0,k}(x_k)}dm_1 & \,\mbox{ if $\,\alpha=\mathbf{0}$}
       			 \vspace{0.01cm} \\
  		       \int_{\mathbb{T}^1} e^{-i\alpha_k\cdot x_k}\ln \frac{1}{p_{0,k}(x_k)}dm_1 & \,\mbox{ if $\,\alpha=(0,\ldots,\alpha_k,\ldots,0)\ne\mathbf{0}$ }
     			 \vspace{0.01cm} \\
  		       0 & \,\mbox{ otherwise }
     			 \vspace{0.01cm} \\
		     \end{array}
		     \right. 
\end{split}		
\end{equation}	 
which is summarized in the following corollary:

\begin{corollary} \label{cor of PNN}
If all components $X_k$ of $X$ are independent 
then there is a unique value of 
\begin{equation*}
	y_{\alpha} = 
		  \left\{
		      \begin{array}{ll}
 		       \sum_{k=1}^d \widehat{\ln \frac{1}{p_{0,k}}}(0) & \,\mbox{ if $\,\alpha=\mathbf{0}$} 
       			 \vspace{0.01cm} \\
  		       \widehat{\ln\frac{1}{p_{0,k}}}(\alpha_k)  & \,\mbox{ if $\,\alpha=(0,\ldots,\alpha_k,\ldots,0)\ne\mathbf{0}$ } 
     			 \vspace{0.01cm} \\
		      0 & \,\mbox{ otherwise}
		     \end{array}
		     \right. 
\end{equation*}
for $1\le k\le d$, which satisfies $(\ref{goal})$, where $p_{0,k}$ is the marginal PDF of $X_k$.
\end{corollary}

\noindent \textit{Note 1.} If every component of $X$ is i.i.d., $p_{0,k}$ is identical, and thus it equals $p_{0,1}$ for every $k$. 
\textit{Note 2.} As a contraposition to the corollary, first, if there exists a multi-index $\alpha$ outside the axis indices such that $y_\alpha \neq 0$, then the components of $X$ are not independent. Second, if there exist two non-trivial multi-indices $\alpha$ and $\alpha'$ in different axes with the same distance from the origin, and $y_\alpha \neq y_{\alpha'}$, it can be concluded that all components of the random variables are not identically distributed.

\smallskip

We would like to mention that by using Legendre polynomials, Chebyshev polynomials, or other orthogonal polynomials instead of $e^{i\alpha\cdot x}$, one might obtain an alternative projection form for $y_\alpha$ when compared to Theorem \ref{PNN by frequency}, leveraging their orthogonality properties (for more information on orthogonal polynomials, refer to \cite{OLBC}).

While Theorem \ref{PNN by frequency} demonstrates the characterization of $y_\alpha$ that allows $p_0$ to be represented as $e^{-E}/Z_E$, $y_\alpha$ is entirely dependent on $p_0$, which still renders Theorem \ref{PNN by frequency} lacking in practicality. In the next section, we will provide a method for inferring $y_\alpha$ from data.

\section{Coefficients from an ergodic process} \label{main section}

As the conclusion section, we will study statistical methods for identifying $y_\alpha$ such that (\ref{goal}) holds based on a sample.
We assume that 
$X=X[n]=(X_1[n],\ldots,X_d[n])$ is an ergodic stochastic process whose realizations reside within the window $\mathbb{T}^d$.

\subsection{System of polynomial series}
We can rewrite $p_0(x)$ as $e^{-E(x,y)}/Z_E$ for a unique $y$ as proven by Theorem \ref{PNN by frequency}.
Since $\E_{p_0}[e^{-i\alpha\cdot X}]=1$ for $\alpha=\mathbf{0}$, let's now consider the case where $\alpha\ne\mathbf{0}$.
It follows that
\begin{equation} \label{base expectation}
\begin{split}
	\E_{p_0}[e^{-i\alpha\cdot X}] 
		&= \int_{\mathbb{T}^d}e^{-i\alpha\cdot x} p_0(x) dx \\
		&= \int_{\mathbb{T}^d}e^{-i\alpha\cdot x} \frac{e^{-E(x,y)}}{Z_E} dx  \\
		&= \frac{1}{Z_E}\sum_{n=0}^\infty\frac{(-1)^n}{n!}\int_{\mathbb{T}^d}e^{-i\alpha\cdot x} E(x,y)^n dx.
\end{split}	
\end{equation}
By rearrangement of series,
\begin{equation*}
	E(x,y)^n = \sum_{\beta^1\in\mathbb{Z}^d}\cdots\sum_{\beta^n\in\mathbb{Z}^d}y_{\beta^1}\cdots y_{\beta^n}e^{i(\beta^1+\cdots+\beta^n)\cdot x},
\end{equation*}
and thus the integration in (\ref{base expectation}) can be calculated as
\begin{equation*}
\begin{split}
	\int_{\mathbb{T}^d}e^{-i\alpha\cdot x}E(x,y)^n dx 
		&= \sum_{\beta^1\in\mathbb{Z}^d}\cdots\sum_{\beta^n\in\mathbb{Z}^d}y_{\beta^1}\cdots y_{\beta^n}
			\int_{\mathbb{T}^d}e^{i(\beta^1+\cdots+\beta^n-\alpha)\cdot x} dx \\
		&= (2\pi)^d\sum_{\substack{\beta^1+\cdots+\beta^n=\alpha \\ \beta^j\in\mathbb{Z}^d } } y_{\beta^1}\cdots y_{\beta^n} \\
		&= (2\pi)^d E(\mathbf{0},y)^{\ast n}(\alpha),
\end{split}		
\end{equation*}
where 
$E(\mathbf{0},y)^{\ast n}(\alpha)$ denotes the $n$-fold iteration of convolution with itself at index $\alpha$, which is known as the convolution power.

For simplicity of notation, if we utilize the notation of the convolution exponential as
\begin{equation} \label{dirac notation}
	\exp^\ast(-E(\mathbf{0},y))(\alpha)=\delta_{\mathbf{0}}(\alpha)+\sum_{n=1}^\infty\frac{(-1)^n}{n!}E(\mathbf{0},y)^{\ast n}(\alpha)
\end{equation}
($\delta_{\mathbf{0}}$ denotes the Dirac delta distribution),
we can express $\E_{p_0}[e^{-i\alpha\cdot X}] $ as follows:
\begin{equation} \label{infinite equations}
	\E_{p_0}[e^{-i\alpha\cdot X}] =  (2\pi)^d\frac{\exp^\ast(-E(\mathbf{0},y))(\alpha) }{Z_E}
				\quad(\alpha\in\mathbb{Z}^d)
\end{equation}
(for more information about the properties of the convolution exponential, refer to \cite{stroock,brouder-patras}).
By the Ergodic theorem, (\ref{infinite equations}) can be expressed as
\begin{equation} \label{homo poly}
	\lim_{M\to\infty}\frac{1}{M}\sum_{n=1}^M e^{-i\alpha\cdot X[n]} 
		=  (2\pi)^d\frac{\exp^\ast(-E(\mathbf{0},y))(\alpha) }{Z_E}
			\quad(\alpha\in\mathbb{Z}^d).
\end{equation}

Especially, if the components of $X$ are independent, then 
\begin{equation*} \label{simple_homo poly}
	p_0(x)=\prod_{k=1}^d p_{0,k}(x_k)=\prod_{k=1}^d p_k(x_k\mid y^{(k)})=p(x\mid y),
\end{equation*}
where $p_{0,k}$, $p_k$ are the marginal PDFs from $p_0$, $p$, respectively, and $y^{(k)}=(y_\eta^{(k)})$ $(\eta\in\mathbb{Z})$.
(Recall that, by Corollary \ref{cor of PNN}, $y_\alpha$ vanishes outside of the axis indices. 
This gives us simpler equations than (\ref{homo poly}):
\begin{equation} \label{homo poly special}
	\lim_{M\to\infty}\frac{1}{M}\sum_{n=1}^M e^{-i \eta X_k[n]} 
		= 2\pi\frac{\exp^\ast(-E(0,y^{(k)}))(\eta) }{Z_{E(0,y^{(k)})}}
			\quad(1\le k\le d,\; \eta\in\mathbb{Z}).
\end{equation}

The left-hand sides of the equations (\ref{homo poly}) and (\ref{homo poly special}) depend only on $X$, while the right-hand sides are comprised of sums of homogeneous polynomials in the components of $y$.
Therefore, the equations reveal themselves as algebraic systems of polynomial series in the components of $y$.
We summarize it as follows, in which the result does not assume a statistical model for the unrevealed PDF.

\begin{theorem} \label{statistical inference} 
Assuming Assumption \ref{A'}, if $X = X[n] = (X_1[n], \ldots, X_d[n]) \in \mathbb{T}^d$ is an ergodic stochastic process, then $(y_\alpha)$, satisfying $(\ref{goal})$, is the unique solution to the system of equations $(\ref{homo poly})$. Moreover, if $X_1, \ldots, X_d$ are independent, then $y$ is the unique solution to the system of equations $(\ref{homo poly special})$.
\end{theorem}

In Theorem \ref{statistical inference}, the equations to be considered are not finite in number; however, the obtained data, in general, consists of only a finite set.
Consequently, we can only anticipate approximate values $y_\alpha$ for a finite set of indices $\alpha$.
For these reasons, it becomes necessary to validate the approximations of the mean and variance of $X$, as well as the shape of $p_0$ based on the values which are  drawn from the hidden $p_0$.

In the following subsection, we will use a random sample drawn from the bivariate normal distribution as a hidden PDF to verify the numerical approximations.

\subsection{Truncated  system of polynomial series}
From a finite amount of data $X[n]$ $(1\le n\le M)$, the number of equations in (\ref{homo poly}) or (\ref{homo poly special}) from Subsection 6.1 will also be truncated. As a result, we will obtain a finite number of approximate $y_\alpha$, denoted as $\hat y_\alpha$, and we will calculate approximate partial sums. Since the convergence of the partial sums of the energy function is guaranteed by Theorem \ref{decaying}, in this subsection, we aim to present truncated equations for this.

Now, let's consider the equation,
\begin{equation} \label{homo poly_pract}
	\frac{1}{M}\sum_{n=1}^M e^{-i\alpha\cdot X[n]} 
		=  \frac{(2\pi)^d}{Z_{E_{N_1}(\mathbf{0},y)}}\left(\delta_{\mathbf{0}}(\alpha)+ \sum_{n=1}^{N_2}\frac{(-1)^n}{n!}E_{N_1}(\mathbf{0},y)^{\ast n}(\alpha)\right)
\end{equation}
for $|\alpha|_\infty\le N_1$, where
$E_{N_1}(\mathbf{0},y)=\sum_{|\alpha|_\infty\le N_1} y_\alpha$ 
represents the partial sum of energy functions at $X=\mathbf{0}$.
If $X_1,X_2,\ldots,N_d$ are independent, then (\ref{homo poly_pract}) is written as
\begin{equation} \label{homo poly_special_pract}
	\frac{1}{M}\sum_{n=1}^M e^{-i \eta X_k[n]} 
		=  \frac{2\pi}{Z_{E_{N_1}(0,y^{(k)})}}\left(\delta_{0}(\eta)+\sum_{n=1}^{N_2}\frac{(-1)^n}{n!}E_{N_1}(0,y^{(k)})^{\ast n}(\eta)\right)
\end{equation}
for $1\le k\le d,\;-N_1\le\eta\le N_1$, where
$E_{N_1}(0,y^{(k)})=\sum_{-N_1\le\eta\le N_1} y^{(k)}_{\eta}$ represents the partial sum at $X_k=0$. (Especially, under the i.i.d. condition, $E_{N_1}(0,y^{(k)})$ is identical for every $k$.)
Here, $M$ denotes the number of realizations $x[n]$ of $X[n]$, and we have the flexibility to choose values for both $N_1$ and $N_2$. The systems described in (\ref{homo poly_pract}) and (\ref{homo poly_special_pract}) consist of $(2N_1+1)^d$-variate and $d(2N_1+1)$-variate polynomials, respectively.

Since $Z_E$ does not depend on $x$, by regarding $y_{\mathbf{0}}$ as $y_{\mathbf{0}}+\ln Z_E$
we solve $(\ref{homo poly_pract})$ or $(\ref{homo poly_special_pract})$ with putting $Z_E=1$. 
The accuracy of solutions $\hat y=(\hat y_\alpha)$ to 
 (\ref{homo poly_pract}) and (\ref{homo poly_special_pract}) will be improved if $M$, $N_1$, and $N_2$ increase according to the Ergodic theorem and Theorem \ref{decaying}.
As a result, $\hat p_0$ serves as a statistical estimator for $p_0$, which is defined as
\begin{equation} \label{PDF estimator}
	\hat p_0(x)= e^{-E(x,\hat y)},
\end{equation}
where $E(x,\hat y)=\sum_{|\alpha|_\infty\le N_1}\hat y_{\alpha}e^{i\alpha\cdot x}$ in which 
$(\hat y_\alpha)$ is regarded as a sequence in $\ell^1$ by trivial extension.

Two systems of (\ref{homo poly_pract}) and (\ref{homo poly_special_pract}) are composed of linear combinations of homogeneous polynomials. 
There have been many research results to approximate their solutions, and there are various algorithms available (\cite{MT, MPR, sturmfels}).
In this experiment, we use the \texttt{solve()} function in MATLAB to find the solution.
Finally, we demonstrate the feasibility of the methods presented with the following numerical results.

\section*{Example}
To avoid the possibility of $p_0(x)=0$ for some $x$, we will consider $p_0(x)+1$ as the target function and adopt (\ref{homo poly_pract}),
since there is no information available on the relation between the components $X_k$ in this example, we will consider (\ref{final approximation system}).
For any $\alpha$, since (\ref{base expectation}) is written as 
\begin{equation} \label{final approximation system}
	\E_{p_0}[e^{-i\alpha\cdot X}] 
	= \int_{\mathbb{T}^d}e^{-i\alpha\cdot x} (p_0(x) + 1)dx -(2\pi)^d\delta_{\mathbf{0}}(\alpha), 
\end{equation}
the solution, namely $(\tilde y_\alpha)$ to 
\begin{equation} \label{final approximation system2}
	\frac{1}{M}\sum_{n=1}^M e^{-i\alpha\cdot X[n]} 
		=  \frac{1}{Z_{E_{N_1}}}\sum_{n=1}^{N_2}\frac{(-1)^n}{n!}E_{N_1}(\mathbf{0},y)^{\ast n}(\alpha)
\end{equation}
is used for approximating $p_0+1$ and then subtract $1$ from the result. 
This technique gives us a simpler implementation than (\ref{homo poly_pract}) in programming and yields slightly more accurate results, along with greater mathematical rigor.
Additionally, we will solve (\ref{final approximation system2}) by inserting $Z_{E_{N_1}}=1$, as mentioned before.
The canonical ensemble form, $e^{-E(x,(\tilde y_\alpha))}$ induced by $(\tilde y_\alpha)$, approximates $p_0+1$, i.e., $\hat p_0:=e^{-E(x,(\tilde y_\alpha))}-1$ is an approximation of  $p_0$. Also,  $(\hat y_\alpha)$ is obtained through the inner product operation of 
$-\ln(e^{-E_{N_1}(x,(\tilde y_\alpha))}-1)$
with 
$e^{i\alpha\cdot x}$. 

We now begin with four sets of random samples, denoted as $(X_1[n], X_2[n])$, each containing counts of $50$, $100$, $200$, and $400$, respectively,
from a bivariate normal distribution $\mathcal{N}_2(\mu, \Sigma)$ with the parameters $\mu = \begin{bmatrix}0 & 0\end{bmatrix}$ and $\Sigma=\begin{bmatrix} 0.25 & 0.2 \\ 0.2 & 0.75\end{bmatrix}$. 
Using the values $N_1 = 5$ and $N_2 = 3$ for the four different cases of $M=50$, $100$, $200$, and $400$,
we obtain the solution set $\tilde y=(\tilde y_\alpha)$, which is depicted in Tables \ref{T1}, \ref{T2}, \ref{T3}, and \ref{T4} for each $M$.
All numbers in each table have been rounded to four decimal places.
However, in the actual calculations, computations were performed up to the default value of sixteen decimal places in the MATLAB program.

\begin{figure}[H]
    \centering 
\begin{subfigure}{0.32\textwidth}
  \includegraphics[trim=32 240 35 220,clip,width=\textwidth]{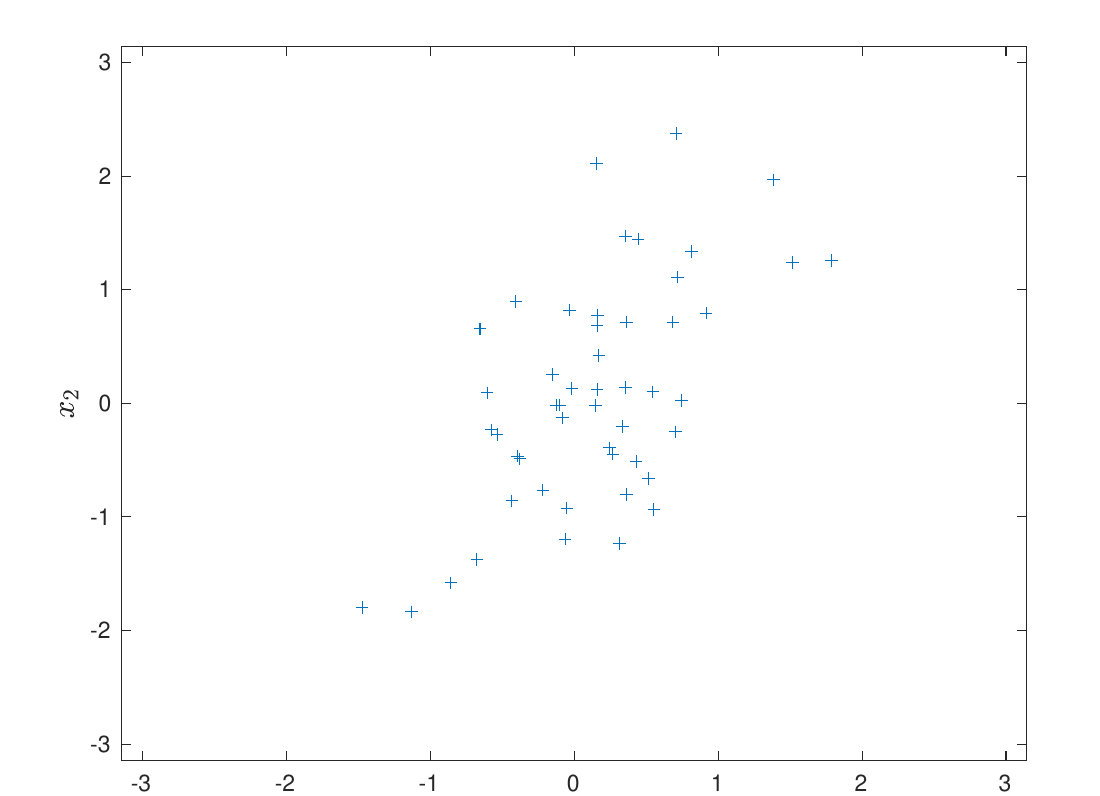}
  \caption{Random sample ($M=50$)}
  \label{fig:1}
\end{subfigure}\hfil 
\begin{subfigure}{0.32\textwidth}
  \includegraphics[trim=32 240 35 220,clip,width=\textwidth]{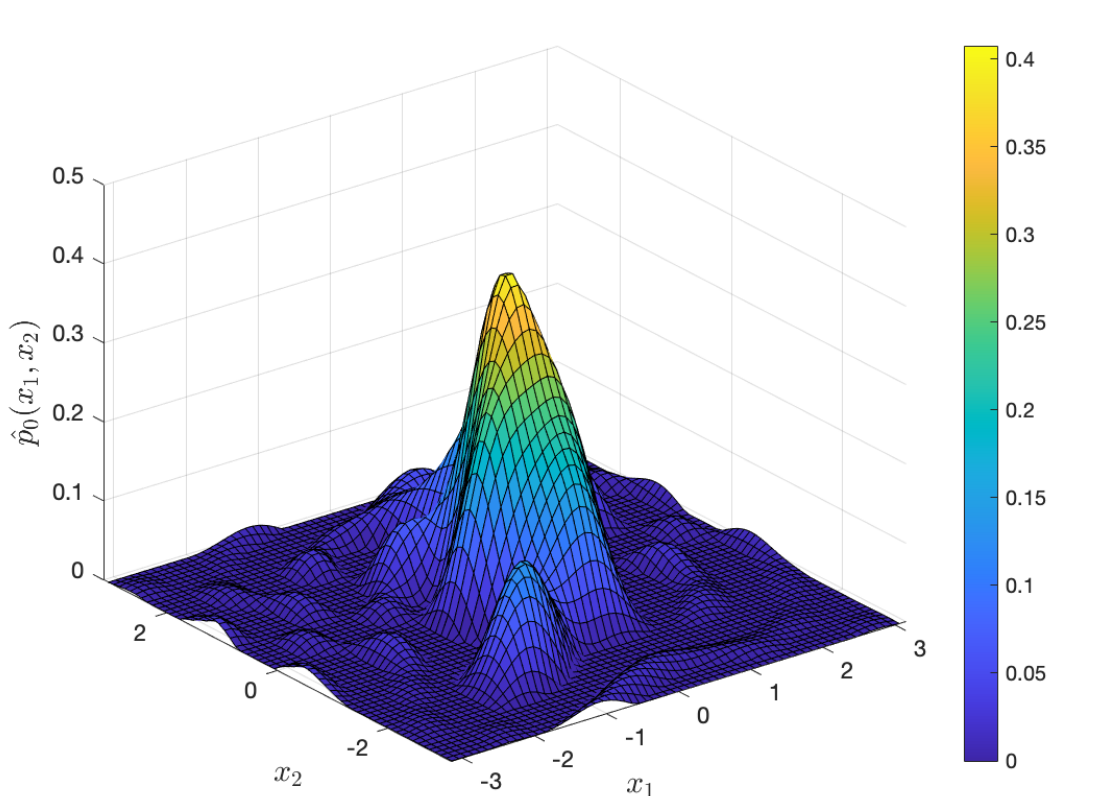}
  \caption{$\hat p_0$ ($M=50$)}
  \label{fig:2}
\end{subfigure}\hfil 
\begin{subfigure}{0.32\textwidth}
  \includegraphics[trim=32 240 35 220,clip,width=\textwidth]{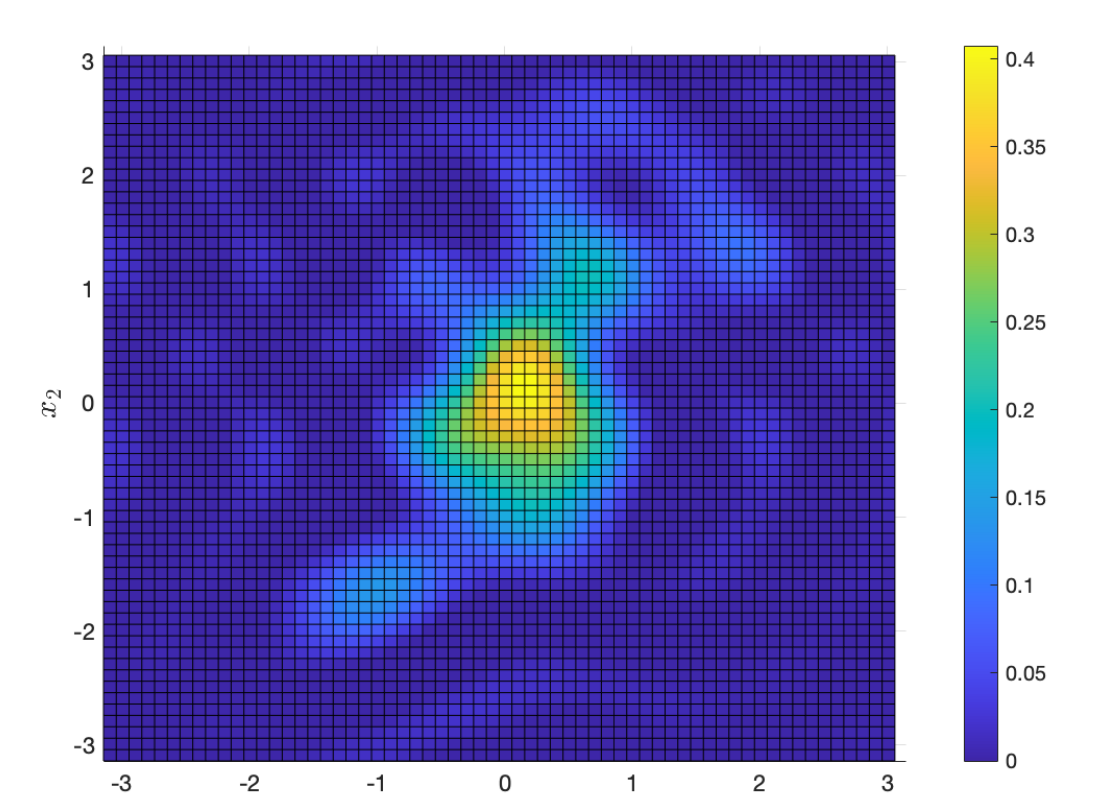}
  \caption{Top view of $\hat p_0$ ($M=50$)}
  \label{fig:3}
\end{subfigure}

\begin{subfigure}{0.32\textwidth}
  \includegraphics[trim=32 240 35 240,clip,width=\textwidth]{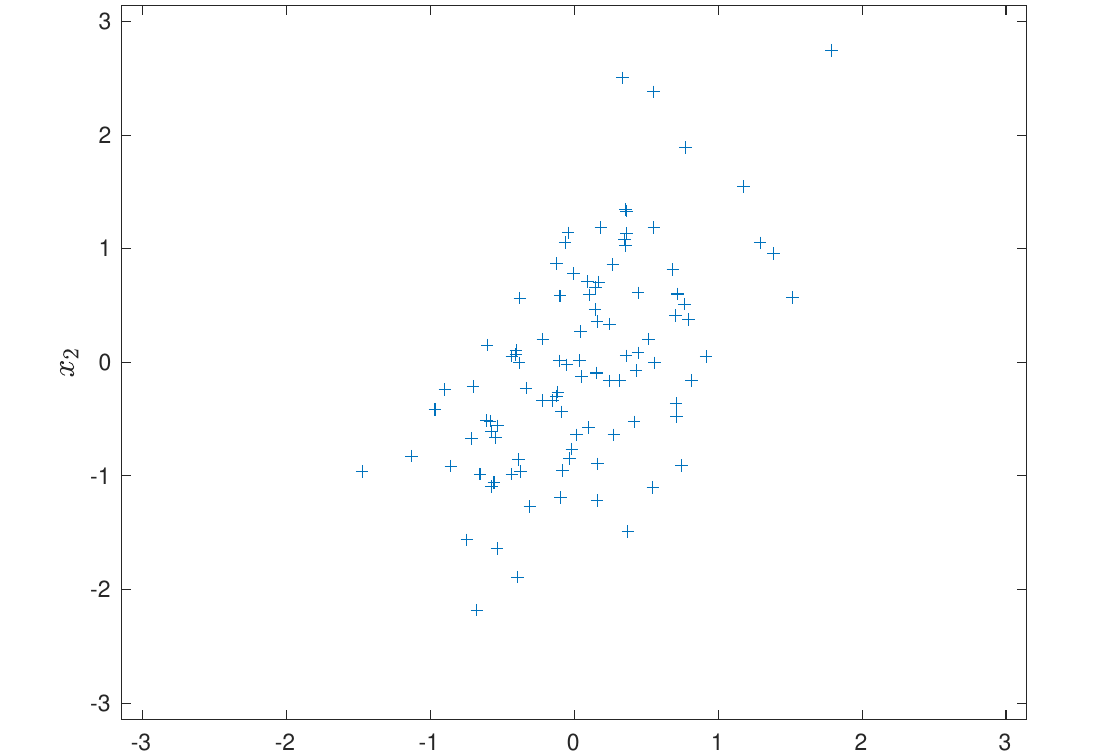}
  \caption{Random sample ($M=100$)}
  \label{fig:4}
\end{subfigure}\hfil 
\begin{subfigure}{0.32\textwidth}
  \includegraphics[trim=32 240 35 240,clip,width=\textwidth]{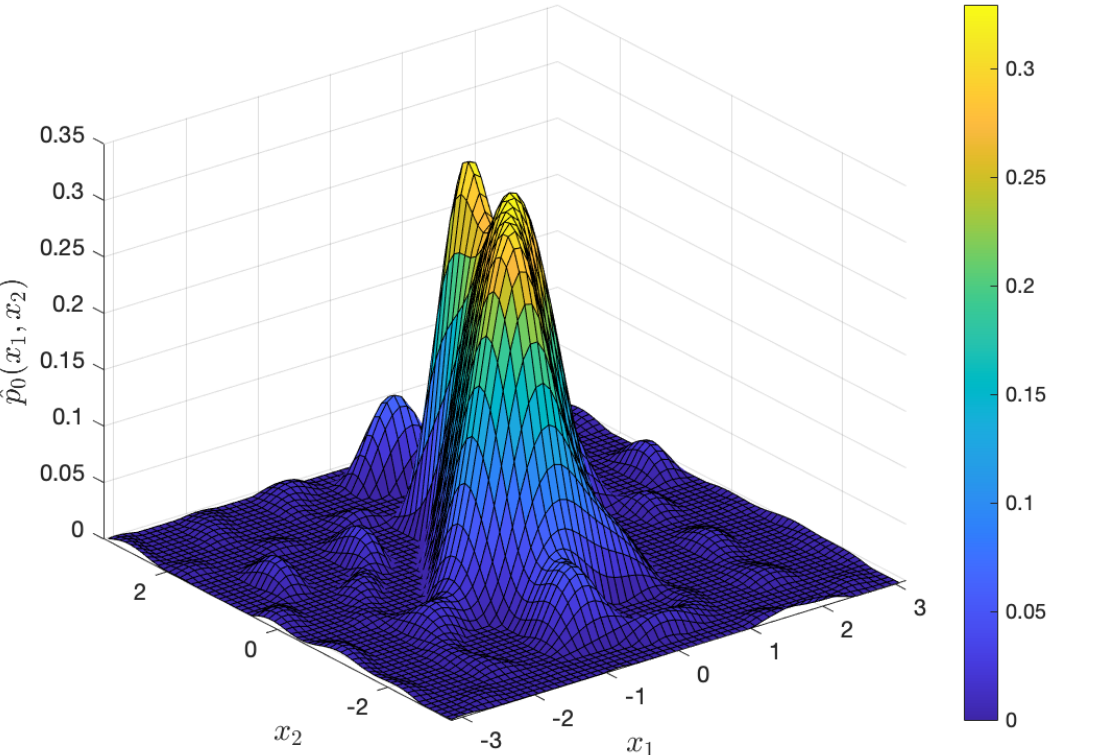}
  \caption{$\hat p_0$ ($M=100$)}
  \label{fig:5}
\end{subfigure}\hfil 
\begin{subfigure}{0.32\textwidth}
  \includegraphics[trim=32 240 35 240,clip,width=\textwidth]{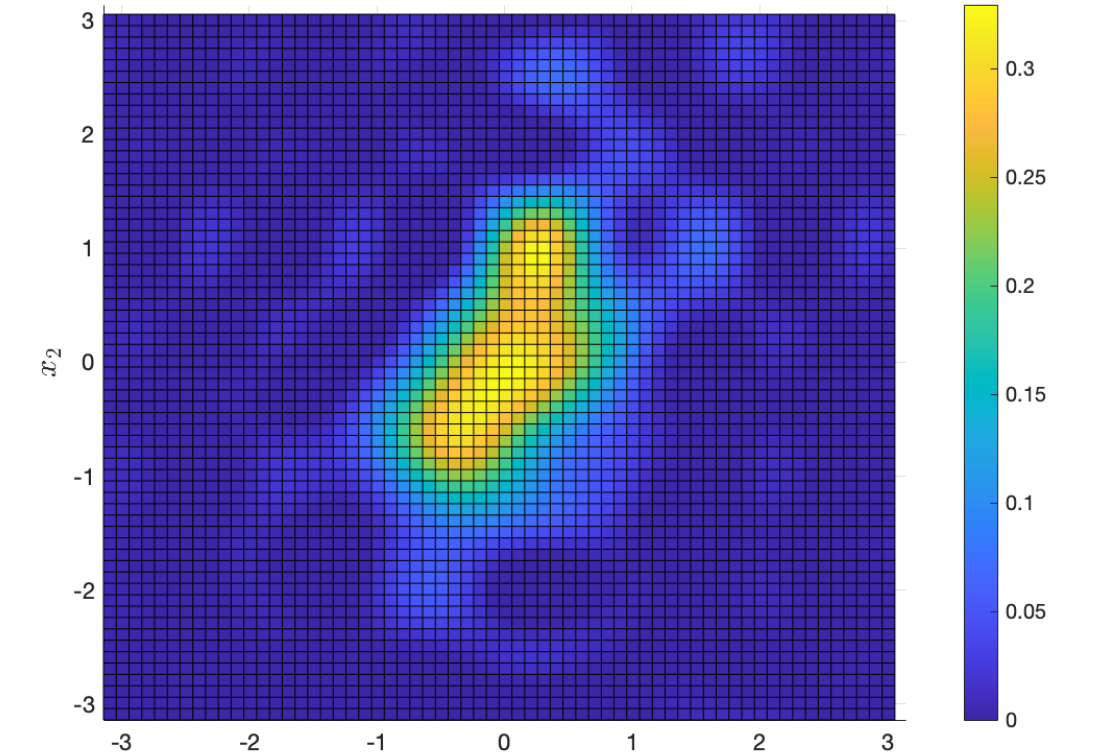}
  \caption{Top view of $\hat p_0$ ($M=100$)}
  \label{fig:6}
\end{subfigure}

\begin{subfigure}{0.32\textwidth}
  \includegraphics[trim=32 240 35 240,clip,width=\textwidth]{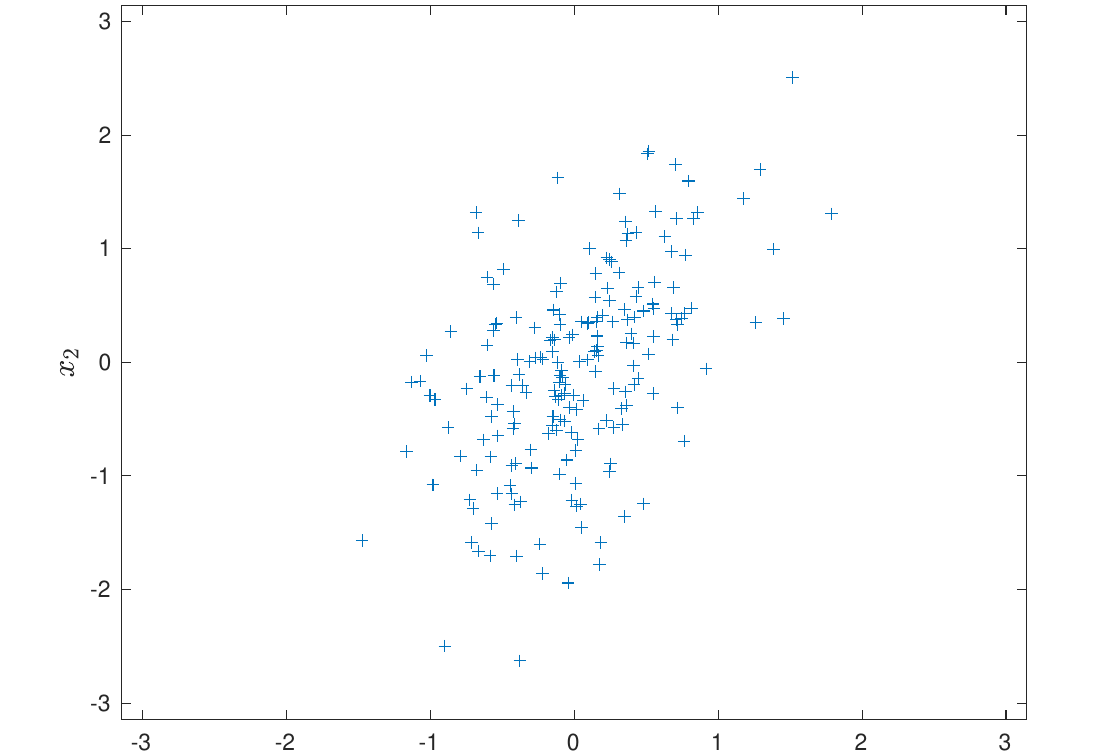}
  \caption{Random sample ($M=200$)}
  \label{fig:7}
\end{subfigure}\hfil 
\begin{subfigure}{0.32\textwidth}
  \includegraphics[trim=32 240 35 240,clip,width=\textwidth]{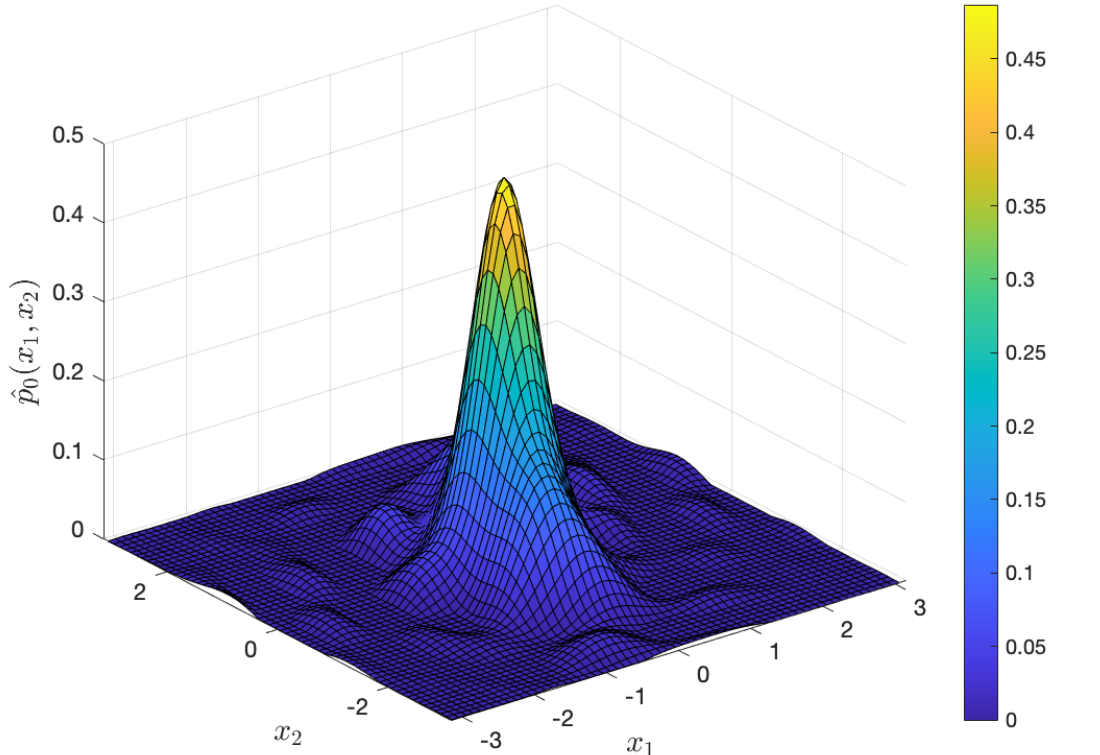}
  \caption{$\hat p_0$ ($M=200$)}
  \label{fig:8}
\end{subfigure}\hfil 
\begin{subfigure}{0.32\textwidth}
  \includegraphics[trim=32 240 35 240,clip,width=\textwidth]{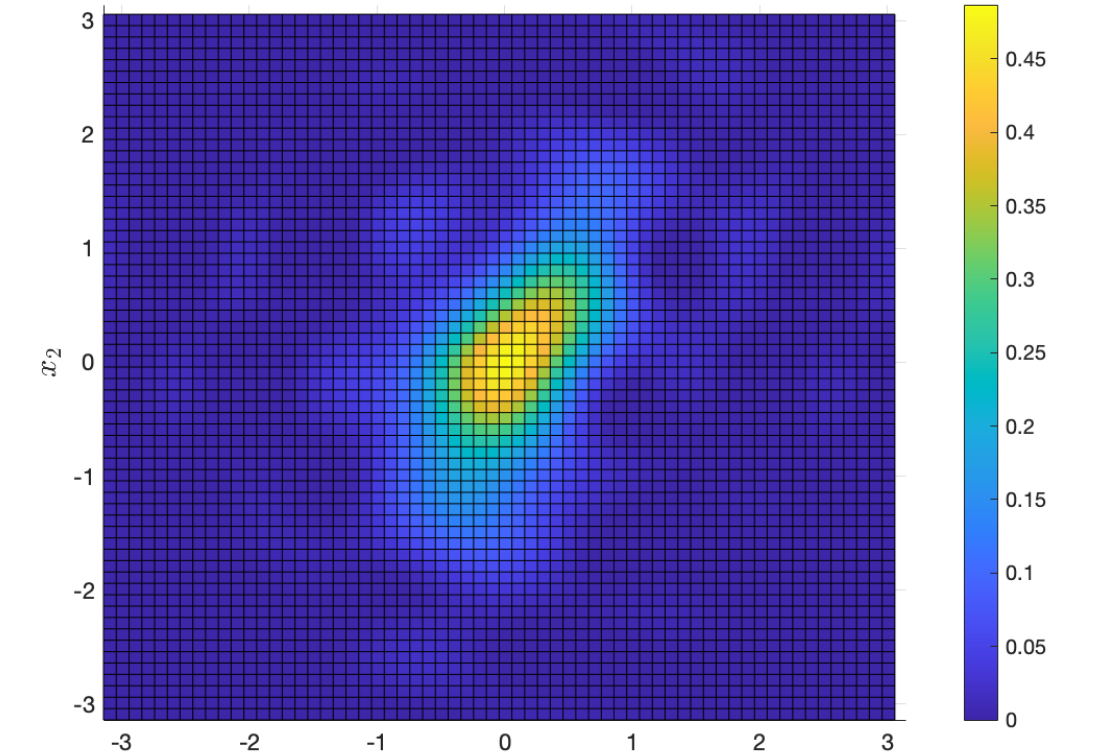}
  \caption{Top view of $\hat p_0$ ($M=200$)}
  \label{fig:9}
\end{subfigure}

\begin{subfigure}{0.32\textwidth}
  \includegraphics[trim=32 240 35 240,clip,width=\textwidth]{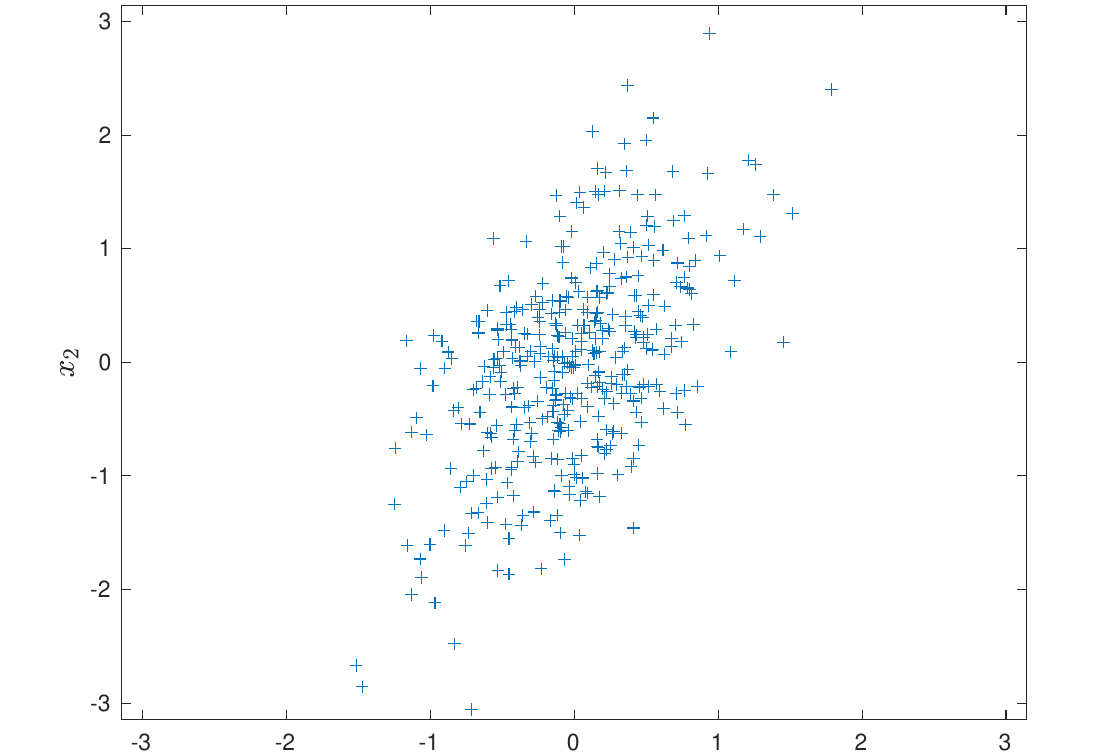}
  \caption{Random sample ($M=400$)}
  \label{fig:10}
\end{subfigure}\hfil 
\begin{subfigure}{0.32\textwidth}
  \includegraphics[trim=32 240 35 240,clip,width=\textwidth]{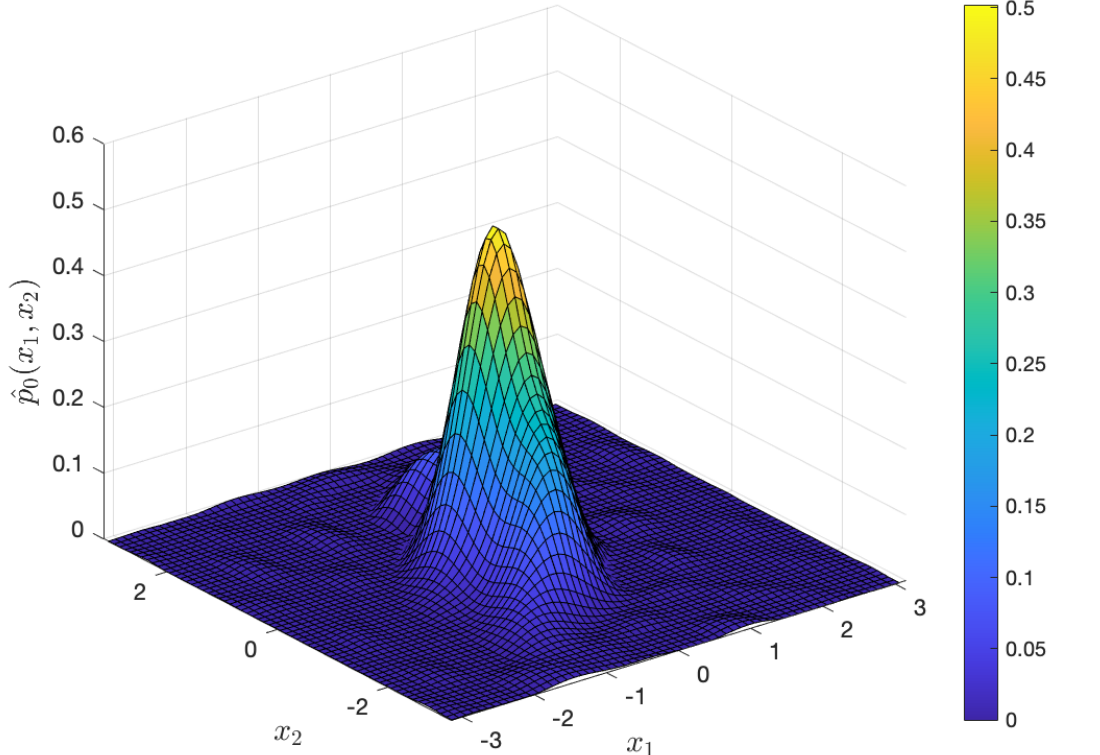}
  \caption{$\hat p_0$ ($M=400$)}
  \label{fig:11}
\end{subfigure}\hfil 
\begin{subfigure}{0.32\textwidth}
  \includegraphics[trim=32 240 35 240,clip,width=\textwidth]{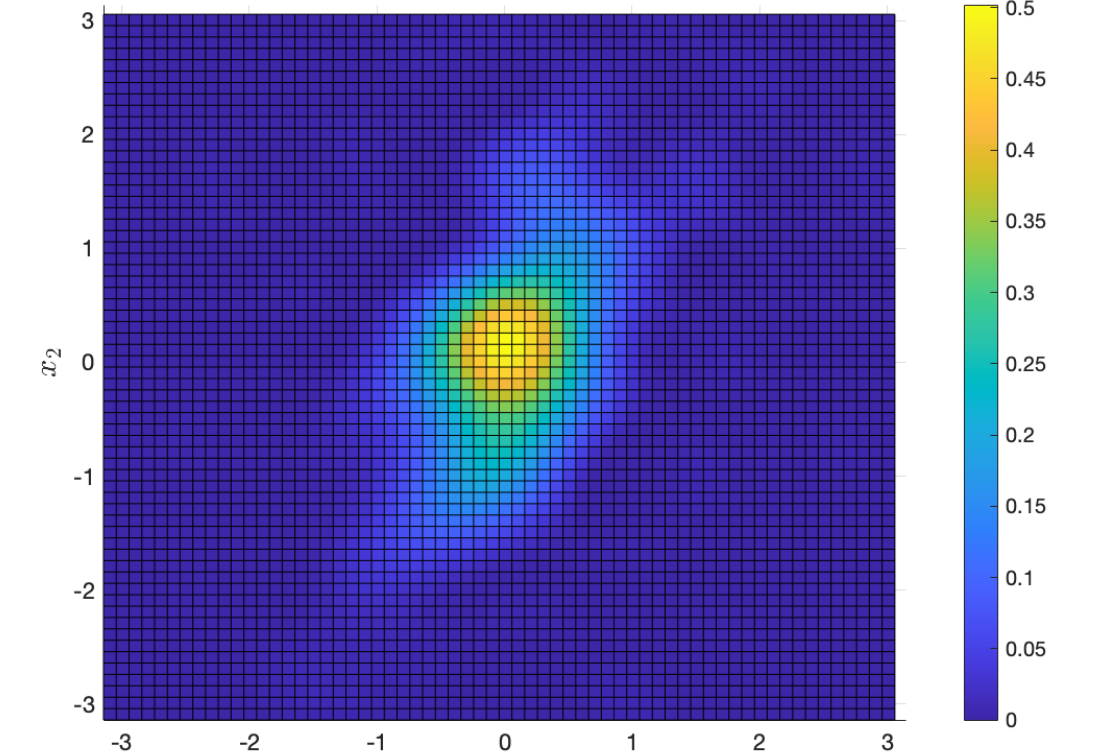}
  \caption{Top view of $\hat p_0$ ($M=400$)}
  \label{fig:12}
\end{subfigure}

\begin{subfigure}{0.32\textwidth}
  \hspace{1cm}
\end{subfigure}\hfil 
\begin{subfigure}{0.32\textwidth}
  \includegraphics[trim=32 240 35 240,clip,width=\textwidth]{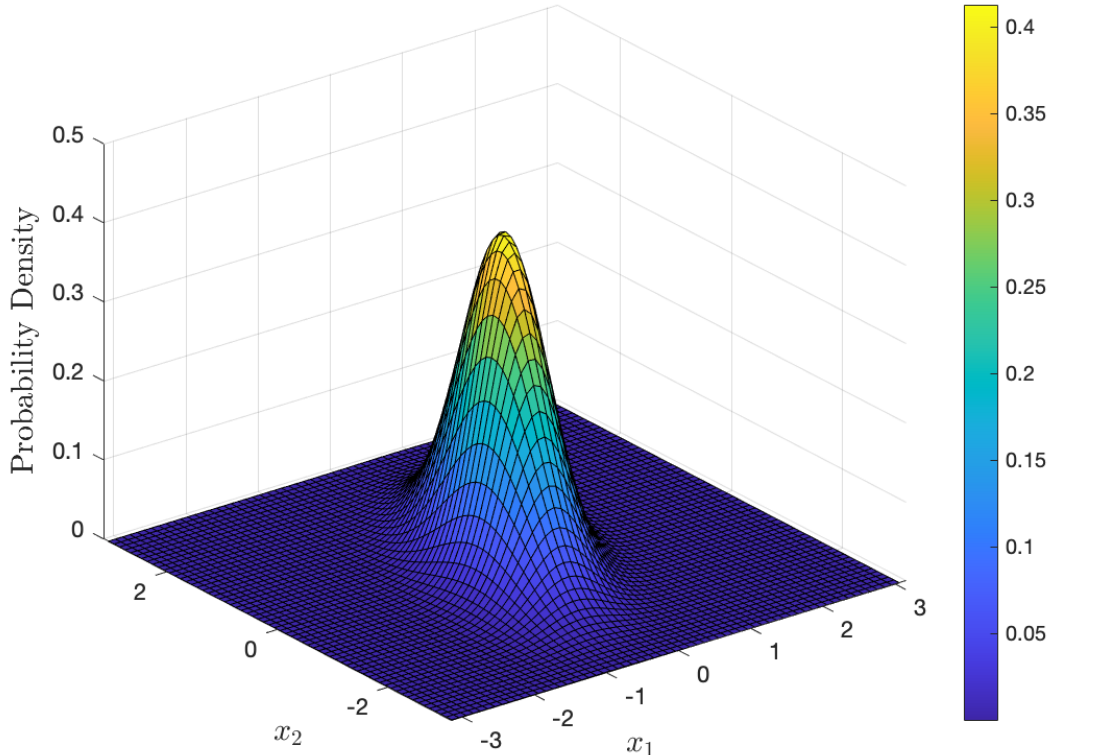}
  \caption{$\mathcal{N}_2(\mu,\Sigma)$}
  \label{fig:13}
\end{subfigure}\hfil 
\begin{subfigure}{0.32\textwidth}
  \includegraphics[trim=32 240 35 240,clip,width=\textwidth]{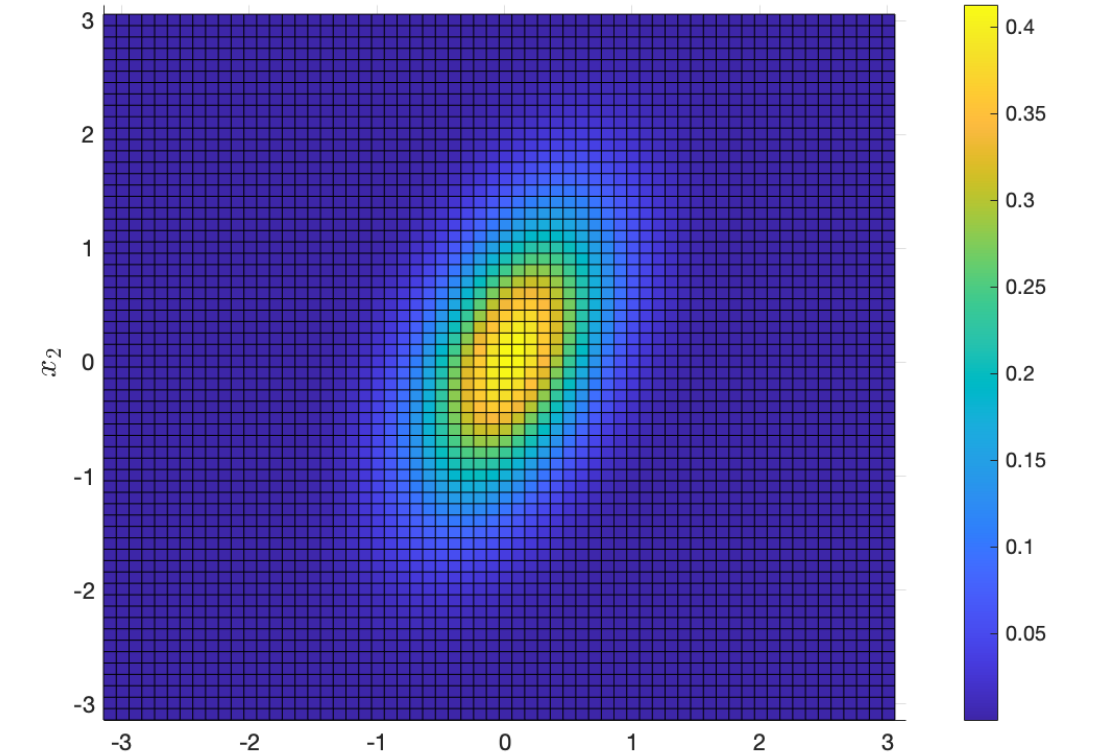}
  \caption{Top view of $\mathcal{N}_2(\mu,\Sigma)$}
  \label{fig:14}
\end{subfigure}
\caption{Approximations of $p_0=\mathcal{N}_2(\mu,\Sigma)$}
\label{fig:images}
\end{figure}

By employing four tables, we computed four approximations of (\ref{PDF estimator}), as illustrated in Figure \ref{fig:images}. As the parameter $M$ increases, the profile of $\hat p_0$ becomes more similar to that of $p_0$. Furthermore, in Table \ref{T5}, we can observe the obtained sample means $\hat\mu$ and sample covariances $\hat\Sigma$, which are calculated from the four $\hat p_0$'s, alongside $\mu$ and $\Sigma$ from $\mathcal{N}_2(\mu, \Sigma)$, the unrevealed PDF. As expected, with the increase in $M$, we can see that both $\hat\mu$ and $\hat\Sigma$ become closer to $\mu$ and $\Sigma$, respectively.

\addtolength{\tabcolsep}{-3.1pt} 
\begin{table}[H]\centering
\ra{1.5}
\begin{tabular}{@{}cccccccccccc|c@{}}\toprule 
$\alpha_1$ & \multicolumn{11}{c@{}}{$\alpha_2$} &  \\
   \cmidrule(l){2-12}
& $-5$ & $-4$ & $-3$ & $-2$ & $-1$ & $0$ & $1$ & $2$ & $3$ & $4$ & \multicolumn{1}{c}{$5$} &  \\ \midrule 
$-5$ &
$\underset{-0.11i}{\overset{0.17}{}}$ &
$\underset{-0.43i}{\overset{0.05}{}}$ &
$\underset{-0.24i}{\overset{-0.22}{}}$ &
$\underset{+0.04i}{\overset{-0.03}{}}$ &
$\underset{-0.04i}{\overset{0.21}{}}$ &
$\underset{-0.27i}{\overset{0.18}{}}$ &
$\underset{-0.11i}{\overset{0.05}{}}$ &
$\underset{+0.22i}{\overset{-0.02}{}}$ &
$\underset{+0.19i}{\overset{-0.22}{}}$ &
$\underset{+0.17i}{\overset{-0.18}{}}$ &
$\underset{+0.14i}{\overset{0.07}{}}$ & \multirow{13}{*}{$\scriptstyle\times 10^{-2}$} \\
$-4$ & 
$\underset{-0.16i}{\overset{0.12}{}}$ &
$\underset{-0.25i}{\overset{-0.10}{}}$ &
$\underset{+0.00i}{\overset{-0.10}{}}$ &
$\underset{+0.11i}{\overset{-0.03}{}}$ &
$\underset{-0.14i}{\overset{-0.11}{}}$ &
$\underset{-0.35i}{\overset{-0.19}{}}$ &
$\underset{-0.08i}{\overset{-0.19}{}}$ &
$\underset{+0.10i}{\overset{-0.14}{}}$ &
$\underset{-0.03i}{\overset{-0.15}{}}$ &
$\underset{+0.08i}{\overset{0.00}{}}$ &
$\underset{+0.16i}{\overset{-0.09}{}}$ & \\
$-3$ & 
$\underset{-0.11i}{\overset{0.01}{}}$ &
$\underset{-0.10i}{\overset{0.01}{}}$ &
$\underset{+0.08i}{\overset{-0.02}{}}$ &
$\underset{+0.11i}{\overset{-0.35}{}}$ &
$\underset{-0.12i}{\overset{-0.48}{}}$ &
$\underset{-0.29i}{\overset{-0.54}{}}$ &
$\underset{-0.13i}{\overset{-0.62}{}}$ &
$\underset{-0.17i}{\overset{-0.40}{}}$ &
$\underset{-0.17i}{\overset{-0.08}{}}$ &
$\underset{+0.20i}{\overset{0.07}{}}$ &
$\underset{+0.19i}{\overset{-0.23}{}}$ & \\
$-2$ & 
$\underset{-0.12i}{\overset{0.01}{}}$ &
$\underset{-0.00i}{\overset{0.10}{}}$ &
$\underset{+0.13i}{\overset{-0.22}{}}$ &
$\underset{+0.19i}{\overset{-0.57}{}}$ &
$\underset{-0.08i}{\overset{-0.60}{}}$ &
$\underset{-0.31i}{\overset{-1.07}{}}$ &
$\underset{-0.24i}{\overset{-1.29}{}}$ &
$\underset{-0.29i}{\overset{-0.58}{}}$ &
$\underset{-0.02i}{\overset{0.04}{}}$ &
$\underset{+0.29i}{\overset{0.09}{}}$ &
$\underset{-0.06i}{\overset{-0.20}{}}$ & \\
$-1$ & 
$\underset{-0.01i}{\overset{-0.04}{}}$ &
$\underset{+0.17i}{\overset{-0.01}{}}$ &
$\underset{+0.19i}{\overset{-0.32}{}}$ &
$\underset{+0.19i}{\overset{-0.39}{}}$ &
$\underset{-0.13i}{\overset{-0.83}{}}$ &
$\underset{-0.27i}{\overset{-1.89}{}}$ &
$\underset{-0.14i}{\overset{-1.68}{}}$ &
$\underset{-0.15i}{\overset{-0.44}{}}$ &
$\underset{+0.09i}{\overset{0.10}{}}$ &
$\underset{+0.06i}{\overset{0.03}{}}$ &
$\underset{-0.37i}{\overset{-0.17}{}}$ & \\
$0$ & 
$\underset{+0.28i}{\overset{-0.13}{}}$ &
$\underset{+0.21i}{\overset{-0.06}{}}$ &
$\underset{+0.05i}{\overset{-0.10}{}}$ &
$\underset{+0.09i}{\overset{-0.25}{}}$ &
$\underset{-0.08i}{\overset{-1.38}{}}$ &
$\underset{+0.00i}{\overset{-2.32}{}}$ &
$\underset{+0.08i}{\overset{-1.38}{}}$ &
$\underset{-0.09i}{\overset{-0.25}{}}$ &
$\underset{-0.05i}{\overset{-0.10}{}}$ &
$\underset{-0.21i}{\overset{-0.06}{}}$ &
$\underset{-0.28i}{\overset{-0.13}{}}$ & \\
$1$ & 
$\underset{+0.37i}{\overset{-0.17}{}}$ &
$\underset{-0.06i}{\overset{0.03}{}}$ &
$\underset{-0.09i}{\overset{0.10}{}}$ &
$\underset{+0.15i}{\overset{-0.44}{}}$ &
$\underset{+0.14i}{\overset{-1.68}{}}$ &
$\underset{+0.27i}{\overset{-1.89}{}}$ &
$\underset{+0.13i}{\overset{-0.83}{}}$ &
$\underset{-0.19i}{\overset{-0.39}{}}$ &
$\underset{-0.19i}{\overset{-0.32}{}}$ &
$\underset{-0.17i}{\overset{-0.01}{}}$ &
$\underset{+0.01i}{\overset{-0.04}{}}$ & \\
$2$ & 
$\underset{+0.06i}{\overset{-0.20}{}}$ &
$\underset{-0.29i}{\overset{0.09}{}}$ &
$\underset{+0.02i}{\overset{0.04}{}}$ &
$\underset{+0.29i}{\overset{-0.58}{}}$ &
$\underset{+0.24i}{\overset{-1.29}{}}$ &
$\underset{+0.31i}{\overset{-1.07}{}}$ &
$\underset{+0.08i}{\overset{-0.60}{}}$ &
$\underset{-0.19i}{\overset{-0.57}{}}$ &
$\underset{-0.13i}{\overset{-0.22}{}}$ &
$\underset{+0.00i}{\overset{0.10}{}}$ &
$\underset{+0.12i}{\overset{0.01}{}}$ & \\
$3$ & 
$\underset{-0.19i}{\overset{-0.23}{}}$ &
$\underset{-0.20i}{\overset{0.07}{}}$ &
$\underset{+0.17i}{\overset{-0.08}{}}$ &
$\underset{+0.17i}{\overset{-0.40}{}}$ &
$\underset{+0.13i}{\overset{-0.62}{}}$ &
$\underset{+0.29i}{\overset{-0.54}{}}$ &
$\underset{+0.12i}{\overset{-0.48}{}}$ &
$\underset{-0.11i}{\overset{-0.35}{}}$ &
$\underset{-0.08i}{\overset{-0.02}{}}$ &
$\underset{+0.10i}{\overset{0.01}{}}$ &
$\underset{+0.11i}{\overset{0.01}{}}$ & \\
$4$ & 
$\underset{-0.16i}{\overset{-0.09}{}}$ &
$\underset{-0.08i}{\overset{0.00}{}}$ &
$\underset{+0.03i}{\overset{-0.15}{}}$ &
$\underset{-0.10i}{\overset{-0.14}{}}$ &
$\underset{+0.08i}{\overset{-0.19}{}}$ &
$\underset{+0.35i}{\overset{-0.19}{}}$ &
$\underset{+0.14i}{\overset{-0.11}{}}$ &
$\underset{-0.11i}{\overset{-0.03}{}}$ &
$\underset{-0.00i}{\overset{-0.10}{}}$ &
$\underset{+0.25i}{\overset{-0.10}{}}$ &
$\underset{+0.16i}{\overset{0.12}{}}$ & \\
$5$ & 
$\underset{-0.14i}{\overset{0.07}{}}$ &
$\underset{-0.17i}{\overset{-0.18}{}}$ &
$\underset{-0.19i}{\overset{-0.22}{}}$ &
$\underset{-0.22i}{\overset{-0.02}{}}$ &
$\underset{+0.11i}{\overset{0.05}{}}$ &
$\underset{+0.27i}{\overset{0.18}{}}$ &
$\underset{+0.04i}{\overset{0.21}{}}$ &
$\underset{-0.04i}{\overset{-0.03}{}}$ &
$\underset{+0.24i}{\overset{-0.22}{}}$ &
$\underset{+0.43i}{\overset{0.05}{}}$ &
$\underset{+0.11i}{\overset{0.17}{}}$ & \\
\bottomrule
\end{tabular}
\caption{$\tilde y_\alpha$ for $M=50$.}
\label{T1}
\end{table}
\addtolength{\tabcolsep}{3.1pt}

\addtolength{\tabcolsep}{-3.1pt} 
\begin{table}[H]\centering
\ra{1.5}
\begin{tabular}{@{}cccccccccccc|c@{}}\toprule
$\alpha_1$ & \multicolumn{11}{c@{}}{$\alpha_2$} &  \\
   \cmidrule(l){2-12}
& $-5$ & $-4$ & $-3$ & $-2$ & $-1$ & $0$ & $1$ & $2$ & $3$ & $4$ & \multicolumn{1}{c}{$5$} &  \\ \midrule 
$-5$ &
$\underset{+0.04i}{\overset{-0.13}{}}$ &
$\underset{+0.11i}{\overset{-0.11}{}}$ &
$\underset{+0.15i}{\overset{0.20}{}}$ &
$\underset{-0.14i}{\overset{0.16}{}}$ &
$\underset{-0.20i}{\overset{0.14}{}}$ &
$\underset{-0.13i}{\overset{0.10}{}}$ &
$\underset{+0.18i}{\overset{-0.28}{}}$ &
$\underset{+0.26i}{\overset{-0.15}{}}$ &
$\underset{+0.14i}{\overset{0.09}{}}$ &
$\underset{+0.07i}{\overset{0.25}{}}$ &
$\underset{-0.19i}{\overset{0.15}{}}$ & \multirow{13}{*}{$\scriptstyle\times 10^{-2}$} \\
$-4$ & 
$\underset{-0.03i}{\overset{-0.15}{}}$ &
$\underset{+0.10i}{\overset{0.01}{}}$ &
$\underset{-0.04i}{\overset{0.28}{}}$ &
$\underset{-0.20i}{\overset{0.07}{}}$ &
$\underset{-0.17i}{\overset{0.10}{}}$ &
$\underset{-0.08i}{\overset{-0.17}{}}$ &
$\underset{+0.18i}{\overset{-0.45}{}}$ &
$\underset{+0.16i}{\overset{-0.20}{}}$ &
$\underset{+0.18i}{\overset{0.06}{}}$ &
$\underset{+0.02i}{\overset{0.15}{}}$ &
$\underset{-0.17i}{\overset{-0.01}{}}$ & \\
$-3$ & 
$\underset{-0.05i}{\overset{-0.20}{}}$ &
$\underset{+0.05i}{\overset{0.09}{}}$ &
$\underset{-0.11i}{\overset{0.18}{}}$ &
$\underset{-0.09i}{\overset{0.03}{}}$ &
$\underset{-0.08i}{\overset{-0.06}{}}$ &
$\underset{-0.03i}{\overset{-0.57}{}}$ &
$\underset{+0.02i}{\overset{-0.79}{}}$ &
$\underset{+0.06i}{\overset{-0.42}{}}$ &
$\underset{+0.11i}{\overset{-0.07}{}}$ &
$\underset{-0.05i}{\overset{-0.02}{}}$ &
$\underset{-0.18i}{\overset{-0.07}{}}$ & \\
$-2$ & 
$\underset{-0.01i}{\overset{-0.22}{}}$ &
$\underset{+0.01i}{\overset{0.09}{}}$ &
$\underset{-0.02i}{\overset{0.15}{}}$ &
$\underset{+0.02i}{\overset{0.04}{}}$ &
$\underset{+0.02i}{\overset{-0.39}{}}$ &
$\underset{-0.08i}{\overset{-1.19}{}}$ &
$\underset{-0.12i}{\overset{-1.35}{}}$ &
$\underset{-0.06i}{\overset{-0.65}{}}$ &
$\underset{-0.01i}{\overset{-0.21}{}}$ &
$\underset{-0.09i}{\overset{-0.03}{}}$ &
$\underset{-0.19i}{\overset{-0.06}{}}$ & \\
$-1$ & 
$\underset{-0.02i}{\overset{-0.17}{}}$ &
$\underset{+0.06i}{\overset{0.12}{}}$ &
$\underset{+0.05i}{\overset{0.21}{}}$ &
$\underset{+0.10i}{\overset{-0.13}{}}$ &
$\underset{+0.09i}{\overset{-0.94}{}}$ &
$\underset{-0.12i}{\overset{-1.95}{}}$ &
$\underset{-0.13i}{\overset{-1.70}{}}$ &
$\underset{-0.17i}{\overset{-0.68}{}}$ &
$\underset{-0.04i}{\overset{-0.15}{}}$ &
$\underset{-0.12i}{\overset{0.12}{}}$ &
$\underset{-0.10i}{\overset{-0.07}{}}$ & \\
$0$ & 
$\underset{-0.01i}{\overset{-0.11}{}}$ &
$\underset{+0.13i}{\overset{0.18}{}}$ &
$\underset{+0.04i}{\overset{0.09}{}}$ &
$\underset{+0.18i}{\overset{-0.45}{}}$ &
$\underset{+0.10i}{\overset{-1.50}{}}$ &
$\underset{+0.00i}{\overset{-2.32}{}}$ &
$\underset{-0.10i}{\overset{-1.50}{}}$ &
$\underset{-0.18i}{\overset{-0.45}{}}$ &
$\underset{-0.04i}{\overset{0.09}{}}$ &
$\underset{-0.13i}{\overset{0.18}{}}$ &
$\underset{+0.01i}{\overset{-0.11}{}}$ & \\
$1$ & 
$\underset{+0.10i}{\overset{-0.07}{}}$ &
$\underset{+0.12i}{\overset{0.12}{}}$ &
$\underset{+0.04i}{\overset{-0.15}{}}$ &
$\underset{+0.17i}{\overset{-0.68}{}}$ &
$\underset{+0.13i}{\overset{-1.70}{}}$ &
$\underset{+0.12i}{\overset{-1.95}{}}$ &
$\underset{-0.09i}{\overset{-0.94}{}}$ &
$\underset{-0.10i}{\overset{-0.13}{}}$ &
$\underset{-0.05i}{\overset{0.21}{}}$ &
$\underset{-0.06i}{\overset{0.12}{}}$ &
$\underset{+0.02i}{\overset{-0.17}{}}$ & \\
$2$ & 
$\underset{+0.19i}{\overset{-0.06}{}}$ &
$\underset{+0.09i}{\overset{-0.03}{}}$ &
$\underset{+0.01i}{\overset{-0.21}{}}$ &
$\underset{+0.06i}{\overset{-0.65}{}}$ &
$\underset{+0.12i}{\overset{-1.35}{}}$ &
$\underset{+0.08i}{\overset{-1.19}{}}$ &
$\underset{-0.02i}{\overset{-0.39}{}}$ &
$\underset{-0.02i}{\overset{0.04}{}}$ &
$\underset{+0.02i}{\overset{0.15}{}}$ &
$\underset{-0.01i}{\overset{0.09}{}}$ &
$\underset{+0.01i}{\overset{-0.22}{}}$ & \\
$3$ & 
$\underset{+0.18i}{\overset{-0.07}{}}$ &
$\underset{+0.05i}{\overset{-0.02}{}}$ &
$\underset{-0.11i}{\overset{-0.07}{}}$ &
$\underset{-0.06i}{\overset{-0.42}{}}$ &
$\underset{-0.02i}{\overset{-0.79}{}}$ &
$\underset{+0.03i}{\overset{-0.57}{}}$ &
$\underset{+0.08i}{\overset{-0.06}{}}$ &
$\underset{+0.09i}{\overset{0.03}{}}$ &
$\underset{+0.11i}{\overset{0.18}{}}$ &
$\underset{-0.05i}{\overset{0.09}{}}$ &
$\underset{+0.05i}{\overset{-0.20}{}}$ & \\
$4$ & 
$\underset{+0.17i}{\overset{-0.01}{}}$ &
$\underset{-0.02i}{\overset{0.15}{}}$ &
$\underset{-0.18i}{\overset{0.06}{}}$ &
$\underset{-0.16i}{\overset{-0.20}{}}$ &
$\underset{-0.18i}{\overset{-0.45}{}}$ &
$\underset{+0.08i}{\overset{-0.17}{}}$ &
$\underset{+0.17i}{\overset{0.10}{}}$ &
$\underset{+0.20i}{\overset{0.07}{}}$ &
$\underset{+0.04i}{\overset{0.28}{}}$ &
$\underset{-0.10i}{\overset{0.01}{}}$ &
$\underset{+0.03i}{\overset{-0.15}{}}$ & \\
$5$ & 
$\underset{+0.19i}{\overset{0.15}{}}$ &
$\underset{-0.07i}{\overset{0.25}{}}$ &
$\underset{-0.14i}{\overset{0.09}{}}$ &
$\underset{-0.26i}{\overset{-0.15}{}}$ &
$\underset{-0.18i}{\overset{-0.28}{}}$ &
$\underset{+0.13i}{\overset{0.10}{}}$ &
$\underset{+0.20i}{\overset{0.14}{}}$ &
$\underset{+0.14i}{\overset{0.16}{}}$ &
$\underset{-0.15i}{\overset{0.20}{}}$ &
$\underset{-0.11i}{\overset{-0.11}{}}$ &
$\underset{-0.04i}{\overset{-0.13}{}}$ & \\
\bottomrule
\end{tabular}
\caption{$\tilde y_\alpha$ for $M=100$.}
\label{T2}
\end{table}
\addtolength{\tabcolsep}{3.1pt} 

\vfill\newpage

\addtolength{\tabcolsep}{-3.1pt} 
\begin{table}[h]\centering
\ra{1.5}
\begin{tabular}{@{}cccccccccccc|c@{}}\toprule
$\alpha_1$ & \multicolumn{11}{c@{}}{$\alpha_2$} &  \\
   \cmidrule(l){2-12}
& $-5$ & $-4$ & $-3$ & $-2$ & $-1$ & $0$ & $1$ & $2$ & $3$ & $4$ & \multicolumn{1}{c}{$5$} &  \\ \midrule 
$-5$ &
$\underset{+0.01i}{\overset{-0.02}{}}$ &
$\underset{+0.13i}{\overset{-0.03}{}}$ &
$\underset{+0.09i}{\overset{-0.07}{}}$ &
$\underset{+0.10i}{\overset{-0.06}{}}$ &
$\underset{+0.07i}{\overset{0.06}{}}$ &
$\underset{-0.08i}{\overset{0.01}{}}$ &
$\underset{-0.18i}{\overset{-0.21}{}}$ &
$\underset{-0.14i}{\overset{-0.33}{}}$ &
$\underset{-0.01i}{\overset{-0.30}{}}$ &
$\underset{+0.03i}{\overset{-0.10}{}}$ &
$\underset{+0.00i}{\overset{0.03}{}}$ & \multirow{13}{*}{$\scriptstyle\times 10^{-2}$} \\
$-4$ & 
$\underset{+0.07i}{\overset{0.01}{}}$ &
$\underset{+0.07i}{\overset{-0.01}{}}$ &
$\underset{-0.03i}{\overset{-0.05}{}}$ &
$\underset{+0.06i}{\overset{-0.00}{}}$ &
$\underset{+0.07i}{\overset{0.01}{}}$ &
$\underset{-0.03i}{\overset{-0.22}{}}$ &
$\underset{-0.18i}{\overset{-0.45}{}}$ &
$\underset{-0.13i}{\overset{-0.43}{}}$ &
$\underset{+0.04i}{\overset{-0.27}{}}$ &
$\underset{+0.01i}{\overset{-0.09}{}}$ &
$\underset{+0.01i}{\overset{0.02}{}}$ & \\
$-3$ & 
$\underset{+0.13i}{\overset{0.08}{}}$ &
$\underset{+0.01i}{\overset{0.03}{}}$ &
$\underset{-0.06i}{\overset{-0.01}{}}$ &
$\underset{+0.06i}{\overset{0.00}{}}$ &
$\underset{+0.13i}{\overset{-0.18}{}}$ &
$\underset{+0.05i}{\overset{-0.63}{}}$ &
$\underset{-0.16i}{\overset{-0.85}{}}$ &
$\underset{-0.07i}{\overset{-0.59}{}}$ &
$\underset{+0.05i}{\overset{-0.28}{}}$ &
$\underset{+0.02i}{\overset{-0.13}{}}$ &
$\underset{+0.07i}{\overset{-0.03}{}}$ & \\
$-2$ & 
$\underset{+0.10i}{\overset{0.15}{}}$ &
$\underset{-0.03i}{\overset{-0.00}{}}$ &
$\underset{-0.03i}{\overset{-0.02}{}}$ &
$\underset{+0.07i}{\overset{-0.09}{}}$ &
$\underset{+0.16i}{\overset{-0.53}{}}$ &
$\underset{+0.03i}{\overset{-1.27}{}}$ &
$\underset{-0.15i}{\overset{-1.38}{}}$ &
$\underset{-0.02i}{\overset{-0.74}{}}$ &
$\underset{+0.05i}{\overset{-0.31}{}}$ &
$\underset{+0.06i}{\overset{-0.18}{}}$ &
$\underset{+0.08i}{\overset{-0.11}{}}$ & \\
$-1$ & 
$\underset{+0.01i}{\overset{0.08}{}}$ &
$\underset{-0.04i}{\overset{-0.11}{}}$ &
$\underset{-0.01i}{\overset{-0.09}{}}$ &
$\underset{+0.05i}{\overset{-0.26}{}}$ &
$\underset{+0.12i}{\overset{-1.04}{}}$ &
$\underset{-0.01i}{\overset{-1.98}{}}$ &
$\underset{-0.10i}{\overset{-1.70}{}}$ &
$\underset{+0.02i}{\overset{-0.72}{}}$ &
$\underset{+0.05i}{\overset{-0.27}{}}$ &
$\underset{+0.07i}{\overset{-0.22}{}}$ &
$\underset{+0.07i}{\overset{-0.17}{}}$ & \\
$0$ & 
$\underset{-0.05i}{\overset{-0.09}{}}$ &
$\underset{-0.05i}{\overset{-0.20}{}}$ &
$\underset{-0.03i}{\overset{-0.19}{}}$ &
$\underset{+0.00i}{\overset{-0.51}{}}$ &
$\underset{+0.07i}{\overset{-1.55}{}}$ &
$\underset{+0.00i}{\overset{-2.30}{}}$ &
$\underset{-0.07i}{\overset{-1.55}{}}$ &
$\underset{-0.00i}{\overset{-0.51}{}}$ &
$\underset{+0.03i}{\overset{-0.19}{}}$ &
$\underset{+0.05i}{\overset{-0.20}{}}$ &
$\underset{+0.05i}{\overset{-0.09}{}}$ & \\
$1$ & 
$\underset{-0.07i}{\overset{-0.17}{}}$ &
$\underset{-0.07i}{\overset{-0.22}{}}$ &
$\underset{-0.05i}{\overset{-0.27}{}}$ &
$\underset{-0.02i}{\overset{-0.72}{}}$ &
$\underset{+0.10i}{\overset{-1.70}{}}$ &
$\underset{+0.01i}{\overset{-1.98}{}}$ &
$\underset{-0.12i}{\overset{-1.04}{}}$ &
$\underset{-0.05i}{\overset{-0.26}{}}$ &
$\underset{+0.01i}{\overset{-0.09}{}}$ &
$\underset{+0.04i}{\overset{-0.11}{}}$ &
$\underset{-0.01i}{\overset{0.08}{}}$ & \\
$2$ & 
$\underset{-0.08i}{\overset{-0.11}{}}$ &
$\underset{-0.06i}{\overset{-0.18}{}}$ &
$\underset{-0.05i}{\overset{-0.31}{}}$ &
$\underset{+0.02i}{\overset{-0.74}{}}$ &
$\underset{+0.15i}{\overset{-1.38}{}}$ &
$\underset{-0.03i}{\overset{-1.27}{}}$ &
$\underset{-0.16i}{\overset{-0.53}{}}$ &
$\underset{-0.07i}{\overset{-0.09}{}}$ &
$\underset{+0.03i}{\overset{-0.02}{}}$ &
$\underset{+0.03i}{\overset{-0.00}{}}$ &
$\underset{-0.10i}{\overset{0.15}{}}$ & \\
$3$ & 
$\underset{-0.07i}{\overset{-0.03}{}}$ &
$\underset{-0.02i}{\overset{-0.13}{}}$ &
$\underset{-0.05i}{\overset{-0.28}{}}$ &
$\underset{+0.07i}{\overset{-0.59}{}}$ &
$\underset{+0.16i}{\overset{-0.85}{}}$ &
$\underset{-0.05i}{\overset{-0.63}{}}$ &
$\underset{-0.13i}{\overset{-0.18}{}}$ &
$\underset{-0.06i}{\overset{0.00}{}}$ &
$\underset{+0.06i}{\overset{-0.01}{}}$ &
$\underset{-0.01i}{\overset{0.03}{}}$ &
$\underset{-0.13i}{\overset{0.08}{}}$ & \\
$4$ & 
$\underset{-0.01i}{\overset{0.02}{}}$ &
$\underset{-0.01i}{\overset{-0.09}{}}$ &
$\underset{-0.04i}{\overset{-0.27}{}}$ &
$\underset{+0.13i}{\overset{-0.43}{}}$ &
$\underset{+0.18i}{\overset{-0.45}{}}$ &
$\underset{+0.03i}{\overset{-0.22}{}}$ &
$\underset{-0.07i}{\overset{0.01}{}}$ &
$\underset{-0.06i}{\overset{-0.00}{}}$ &
$\underset{+0.03i}{\overset{-0.05}{}}$ &
$\underset{-0.07i}{\overset{-0.01}{}}$ &
$\underset{-0.07i}{\overset{0.01}{}}$ & \\
$5$ & 
$\underset{-0.00i}{\overset{0.03}{}}$ &
$\underset{-0.03i}{\overset{-0.10}{}}$ &
$\underset{+0.01i}{\overset{-0.30}{}}$ &
$\underset{+0.14i}{\overset{-0.33}{}}$ &
$\underset{+0.18i}{\overset{-0.21}{}}$ &
$\underset{+0.08i}{\overset{0.01}{}}$ &
$\underset{-0.07i}{\overset{0.06}{}}$ &
$\underset{-0.10i}{\overset{-0.06}{}}$ &
$\underset{-0.09i}{\overset{-0.07}{}}$ &
$\underset{-0.13i}{\overset{-0.03}{}}$ &
$\underset{-0.01i}{\overset{-0.02}{}}$ & \\
\bottomrule
\end{tabular}
\caption{$\tilde y_\alpha$ for $M=200$.}
\label{T3}
\end{table}
\addtolength{\tabcolsep}{3.1pt}

\addtolength{\tabcolsep}{-3.1pt} 
\begin{table}[H]\centering
\ra{1.5}
\begin{tabular}{@{}cccccccccccc|c@{}}\toprule
$\alpha_1$ & \multicolumn{11}{c@{}}{$\alpha_2$} &  \\
   \cmidrule(l){2-12}
& $-5$ & $-4$ & $-3$ & $-2$ & $-1$ & $0$ & $1$ & $2$ & $3$ & $4$ & \multicolumn{1}{c}{$5$} &  \\ \midrule 
$-5$ &
$\underset{-0.08i}{\overset{0.05}{}}$ &
$\underset{-0.03i}{\overset{-0.06}{}}$ &
$\underset{+0.01i}{\overset{-0.06}{}}$ &
$\underset{+0.08i}{\overset{0.00}{}}$ &
$\underset{+0.02i}{\overset{-0.02}{}}$ &
$\underset{-0.07i}{\overset{-0.13}{}}$ &
$\underset{-0.03i}{\overset{-0.17}{}}$ &
$\underset{+0.02i}{\overset{-0.02}{}}$ &
$\underset{-0.09i}{\overset{0.04}{}}$ &
$\underset{-0.00i}{\overset{0.03}{}}$ &
$\underset{+0.06i}{\overset{0.02}{}}$ & \multirow{13}{*}{$\scriptstyle\times 10^{-2}$} \\
$-4$ & 
$\underset{-0.06i}{\overset{-0.05}{}}$ &
$\underset{-0.03i}{\overset{-0.08}{}}$ &
$\underset{+0.02i}{\overset{-0.06}{}}$ &
$\underset{+0.08i}{\overset{-0.03}{}}$ &
$\underset{+0.02i}{\overset{-0.12}{}}$ &
$\underset{-0.04i}{\overset{-0.33}{}}$ &
$\underset{+0.02i}{\overset{-0.43}{}}$ &
$\underset{-0.04i}{\overset{-0.23}{}}$ &
$\underset{-0.07i}{\overset{-0.02}{}}$ &
$\underset{+0.11i}{\overset{0.04}{}}$ &
$\underset{+0.13i}{\overset{0.03}{}}$ & \\
$-3$ & 
$\underset{-0.03i}{\overset{-0.09}{}}$ &
$\underset{-0.02i}{\overset{-0.08}{}}$ &
$\underset{+0.03i}{\overset{-0.11}{}}$ &
$\underset{+0.09i}{\overset{-0.10}{}}$ &
$\underset{+0.04i}{\overset{-0.32}{}}$ &
$\underset{+0.02i}{\overset{-0.71}{}}$ &
$\underset{+0.03i}{\overset{-0.90}{}}$ &
$\underset{-0.06i}{\overset{-0.50}{}}$ &
$\underset{+0.03i}{\overset{-0.08}{}}$ &
$\underset{+0.18i}{\overset{0.03}{}}$ &
$\underset{+0.17i}{\overset{0.03}{}}$ & \\
$-2$ & 
$\underset{-0.05i}{\overset{-0.07}{}}$ &
$\underset{-0.02i}{\overset{-0.15}{}}$ &
$\underset{+0.04i}{\overset{-0.17}{}}$ &
$\underset{+0.10i}{\overset{-0.24}{}}$ &
$\underset{+0.06i}{\overset{-0.65}{}}$ &
$\underset{+0.05i}{\overset{-1.33}{}}$ &
$\underset{-0.00i}{\overset{-1.44}{}}$ &
$\underset{-0.03i}{\overset{-0.69}{}}$ &
$\underset{+0.10i}{\overset{-0.15}{}}$ &
$\underset{+0.18i}{\overset{-0.03}{}}$ &
$\underset{+0.15i}{\overset{0.01}{}}$ & \\
$-1$ & 
$\underset{-0.07i}{\overset{-0.07}{}}$ &
$\underset{-0.04i}{\overset{-0.21}{}}$ &
$\underset{+0.01i}{\overset{-0.21}{}}$ &
$\underset{+0.08i}{\overset{-0.43}{}}$ &
$\underset{+0.05i}{\overset{-1.12}{}}$ &
$\underset{+0.03i}{\overset{-2.00}{}}$ &
$\underset{-0.02i}{\overset{-1.72}{}}$ &
$\underset{+0.01i}{\overset{-0.72}{}}$ &
$\underset{+0.08i}{\overset{-0.21}{}}$ &
$\underset{+0.14i}{\overset{-0.11}{}}$ &
$\underset{+0.10i}{\overset{-0.05}{}}$ & \\
$0$ & 
$\underset{-0.07i}{\overset{-0.08}{}}$ &
$\underset{-0.09i}{\overset{-0.19}{}}$ &
$\underset{-0.03i}{\overset{-0.23}{}}$ &
$\underset{+0.02i}{\overset{-0.61}{}}$ &
$\underset{+0.03i}{\overset{-1.56}{}}$ &
$\underset{+0.00i}{\overset{-2.30}{}}$ &
$\underset{-0.03i}{\overset{-1.56}{}}$ &
$\underset{-0.02i}{\overset{-0.61}{}}$ &
$\underset{+0.03i}{\overset{-0.23}{}}$ &
$\underset{+0.09i}{\overset{-0.19}{}}$ &
$\underset{+0.07i}{\overset{-0.08}{}}$ & \\
$1$ & 
$\underset{-0.10i}{\overset{-0.05}{}}$ &
$\underset{-0.14i}{\overset{-0.11}{}}$ &
$\underset{-0.08i}{\overset{-0.21}{}}$ &
$\underset{-0.01i}{\overset{-0.72}{}}$ &
$\underset{+0.02i}{\overset{-1.72}{}}$ &
$\underset{-0.03i}{\overset{-2.00}{}}$ &
$\underset{-0.05i}{\overset{-1.12}{}}$ &
$\underset{-0.08i}{\overset{-0.43}{}}$ &
$\underset{-0.01i}{\overset{-0.21}{}}$ &
$\underset{+0.04i}{\overset{-0.21}{}}$ &
$\underset{+0.07i}{\overset{-0.07}{}}$ & \\
$2$ & 
$\underset{-0.15i}{\overset{0.01}{}}$ &
$\underset{-0.18i}{\overset{-0.03}{}}$ &
$\underset{-0.10i}{\overset{-0.15}{}}$ &
$\underset{+0.03i}{\overset{-0.69}{}}$ &
$\underset{+0.00i}{\overset{-1.44}{}}$ &
$\underset{-0.05i}{\overset{-1.33}{}}$ &
$\underset{-0.06i}{\overset{-0.65}{}}$ &
$\underset{-0.10i}{\overset{-0.24}{}}$ &
$\underset{-0.04i}{\overset{-0.17}{}}$ &
$\underset{+0.02i}{\overset{-0.15}{}}$ &
$\underset{+0.05i}{\overset{-0.07}{}}$ & \\
$3$ & 
$\underset{-0.17i}{\overset{0.03}{}}$ &
$\underset{-0.18i}{\overset{0.03}{}}$ &
$\underset{-0.03i}{\overset{-0.08}{}}$ &
$\underset{+0.06i}{\overset{-0.50}{}}$ &
$\underset{-0.03i}{\overset{-0.90}{}}$ &
$\underset{-0.02i}{\overset{-0.71}{}}$ &
$\underset{-0.04i}{\overset{-0.32}{}}$ &
$\underset{-0.09i}{\overset{-0.10}{}}$ &
$\underset{-0.03i}{\overset{-0.11}{}}$ &
$\underset{+0.02i}{\overset{-0.08}{}}$ &
$\underset{+0.03i}{\overset{-0.09}{}}$ & \\
$4$ & 
$\underset{-0.13i}{\overset{0.03}{}}$ &
$\underset{-0.11i}{\overset{0.04}{}}$ &
$\underset{+0.07i}{\overset{-0.02}{}}$ &
$\underset{+0.04i}{\overset{-0.23}{}}$ &
$\underset{-0.02i}{\overset{-0.43}{}}$ &
$\underset{+0.04i}{\overset{-0.33}{}}$ &
$\underset{-0.02i}{\overset{-0.12}{}}$ &
$\underset{-0.08i}{\overset{-0.03}{}}$ &
$\underset{-0.02i}{\overset{-0.06}{}}$ &
$\underset{+0.03i}{\overset{-0.08}{}}$ &
$\underset{+0.06i}{\overset{-0.05}{}}$ & \\
$5$ & 
$\underset{-0.06i}{\overset{0.02}{}}$ &
$\underset{+0.00i}{\overset{0.03}{}}$ &
$\underset{+0.09i}{\overset{0.04}{}}$ &
$\underset{-0.02i}{\overset{-0.02}{}}$ &
$\underset{+0.03i}{\overset{-0.17}{}}$ &
$\underset{+0.07i}{\overset{-0.13}{}}$ &
$\underset{-0.02i}{\overset{-0.02}{}}$ &
$\underset{-0.08i}{\overset{0.00}{}}$ &
$\underset{-0.01i}{\overset{-0.06}{}}$ &
$\underset{+0.03i}{\overset{-0.06}{}}$ &
$\underset{+0.08i}{\overset{0.05}{}}$ & \\
\bottomrule
\end{tabular}
\caption{$\tilde y_\alpha$ for $M=400$.}
\label{T4}
\end{table}
\addtolength{\tabcolsep}{3.1pt} 

\vfill\newpage

\begin{table}[H]\centering
\ra{1.20}
\begin{tabular}{lccccccl}\toprule
& & \multicolumn{2}{c}{$\hat\mu$} & & \multicolumn{3}{c}{$\hat\Sigma$}
\\\cmidrule(lr){3-4}\cmidrule(lr){6-8}
          & & $\hat\mu_{1}$  & $\hat\mu_{2}$    & & $\hat\Sigma_{11}$  & $\hat\Sigma_{12}$ & $\hat\Sigma_{22}$ \\\midrule
$M=50$  &  & $0.1531$ & $0.0787$  & & 0.3968 & 0.3768 & 0.9604 \\
$M=100$ & & $0.0706$ & $0.0067$  & & 0.3119 & 0.3119 & 0.8536 \\
$M=200$ & & $0.0186$ & $0.0539$  & & 0.2972 & 0.2654 & 0.7604 \\
$M=400$  & & $-0.0115$ & $-0.0171$  & & 0.2659 & 0.2437 & 0.7591 \\\midrule 
& & \multicolumn{2}{c}{$\mu$} & & \multicolumn{3}{c}{$\Sigma$}
\\\cmidrule(lr){3-4}\cmidrule(lr){6-8}
          & & $\mu_{1}$  & $\mu_{2}$    & & $\Sigma_{11}$  & $\Sigma_{12}$ & $\Sigma_{22}$ \\\midrule
$\mathcal{N}_2(\mu,\Sigma)$  & & $0$ & $0$  & & $0.25$ & $0.2$ & $0.75$  \\\bottomrule
\end{tabular}
\caption{Estimations $\hat\mu$ and $\hat\Sigma$ of $\mu$ and $\Sigma$, resp.}
\label{T5}
\end{table}

\section*{Conclusion}
In this mathematical study, we delve into the realm of statistical inference and introduce a novel approach to variational non-Bayesian inference. Most significantly, we propose a new method for uniquely determining the hidden PDF solely from a random sample while leveraging theory.

Beyond merely approximating the shape of the unrevealed PDF, we yield results that provide practical assistance in inferring important moments, such as the mean and variance, from the estimated PDF.

These days, methods utilizing artificial intelligence are widely employed, but they generally exhibit a drawback of depending on initial conditions, i.e., a prior distribution, and iterative computations, i.e., a backpropagation. The limitations of such approaches lie in their reliance on these conditions to achieve global optimization. We emphasize the significant contribution of the proposed method to predictive and classification models, even without any information on populations, highlighting its potential applications in various domains such as finance, economics, weather forecasting, and machine learning, all of which present unique challenges.

The most well-known method for approximating hidden PDFs is the variational approach  in a Bayesian context.
In contrast, our study extends this problem by seeking to precisely determine unknown PDFs through a system of equations.
Our approach includes proving the Fréchet differentiability of entropy to establish the uniqueness of the energy function space in the Wiener algebra. We then derive the unique determination of the energy function through the minimization of KL-divergence. 

Leveraging the Ergodic theorem, we elucidate that solutions to equations comprising polynomial function series are the coefficients of the energy function and numerically substantiate the convergence of partial sums of the energy function obtained from a finite number of equations.

In summary, our mathematical exploration has unveiled the potential of variational non-Bayesian inference in Wiener space. We anticipate that these mathematical ideas can offer an innovative framework for probability density estimation and predictive modeling. Ultimately, our research emphasizes its potential contribution to reshaping mathematical horizons and expanding the boundaries of knowledge in the field of statistical methodology.

\bibliography{my_bibtex}
\bibliographystyle{amsplain}

\end{document}